\documentclass[11pt]{article}

\usepackage{fix-cm}
\usepackage{amsmath,amssymb,amsthm}
\usepackage{newtxtext}
\usepackage{newtxmath}

\usepackage{bm}
\usepackage[margin=1in]{geometry}
\usepackage{setspace}
\usepackage{graphicx}
\usepackage{xcolor}
\usepackage{url}
\makeatletter
\g@addto@macro\UrlBreaks{\do\a\do\b\do\c\do\d\do\e\do\f\do\g\do\h\do\i\do\j\do\k\do\l\do\m\do\n\do\o\do\p\do\q\do\r\do\s\do\t\do\u\do\v\do\w\do\x\do\y\do\z}
\makeatother
\usepackage[round,authoryear]{natbib}
\bibpunct[, ]{(}{)}{,}{a}{}{,}%
\usepackage{algorithm}
\usepackage{algorithmic}
\usepackage{booktabs}
\usepackage{multirow}
\usepackage{tikz}
\usetikzlibrary{arrows.meta,positioning,shapes.geometric}
\usepackage{empheq}
\usepackage[colorlinks=true,linkcolor=blue!60!black,citecolor=blue!60!black,urlcolor=blue!60!black]{hyperref}

\theoremstyle{plain}
\newtheorem{theorem}{Theorem}
\newtheorem{proposition}{Proposition}
\newtheorem{corollary}{Corollary}

\DeclareMathOperator*{\argmax}{arg\,max}

\newcommand{\cA}{{\mathcal{A}}}
\newcommand{\cB}{{\mathcal{B}}}

\newcommand{\cG}{{\mathcal{G}}}

\newcommand{\cD}{{\mathcal{D}}}

\newcommand{\cL}{{\mathcal{L}}}

\newcommand{\cK}{{\mathcal{K}}}
\newcommand{\cF}{{\mathcal{F}}}
\newcommand{\cN}{{\mathcal{N}}}

\newcommand{\cI}{{\mathcal{I}}}

\newcommand{\by}{\textbf{y}}

\newcommand{\bq}{\textbf{q}}
\newcommand{\bu}{\textbf{u}}

\newcommand{\ba}{\textbf{a}}

\newcommand{\bz}{\textbf{z}}

\newcommand{\bh}{\textbf{h}}
\newcommand{\bk}{\textbf{k}}

\newcommand{\bt}{\textbf{t}}
 
\newcommand{\bbt}{\pmb{\beta}}

\newcommand{\bmu}{{\pmb{\mu}}}

\newcommand{\bgamma}{\pmb{\gamma}}
\newcommand{\blambda}{\pmb{\lambda}}
\newcommand{\brho}{\pmb{\rho}}

\newcommand{\NL}{\mathrm{NL}}

\begin{document}

\title{Globally Optimal Estimation of Nested Logit Models}
\author{%
Hoang Giang Pham\thanks{Singapore Management University, Singapore.}
\and
Tien Mai\thanks{Singapore Management University, Singapore. Email: \texttt{atmai@smu.edu.sg}}
}
\date{}

\maketitle

\begin{abstract}
Nested logit is widely used in operational demand models, yet jointly estimating its utility coefficients and nest dissimilarity parameters requires nonconvex optimization. Standard estimation routines offer no certificate of global optimality, leaving a gap between globally optimizing decisions under nested-logit demand and reliably estimating that demand model. We develop, to our knowledge, the first globally convergent optimization framework with an explicit global-accuracy guarantee for joint nested-logit estimation. Our central complexity result confines all exponential dependence to the number of free nest dissimilarity parameters: computation is polynomial in the numbers of observations, alternatives, and utility coefficients. Consequently, under the stated compact-domain and relaxation-error conditions, fixing the number of nests yields an additive polynomial-time approximation scheme (PTAS). The framework exploits an exact exponential-cone formulation for fixed dissimilarities and a spatial branch-and-bound algorithm that searches only over the free dissimilarities. Valid likelihood bounds provide an anytime certificate and finite termination at any prescribed additive average-log-likelihood tolerance. Across eight benchmark datasets and 240 estimation instances, the method certifies every reported solution and exposes inferior local estimates despite normal solver termination. These results establish a computational foundation for globally certified nested-logit estimation.
\end{abstract}


\section{Introduction}

Choice models are embedded in a broad range of operations research (OR) decisions. In revenue management, demand estimates determine which products to offer and which capacity to reserve \citep{Talluri2004revenue,rusmevichientong2010dynamic}. In pricing, they determine markups and product-line design \citep{wang2012capacitated,chen2020capacitated}. In transportation and facility location, they determine predicted flows and captured demand \citep{McFadden1974,Train2009,FLO_Haase2014comparison,Dam2021submodularity}. Across these applications, the estimated substitution pattern is not merely descriptive: it is an input to an optimization model and can change the recommended decision.

The multinomial logit (MNL) model is especially attractive because its probabilities are closed form and its maximum-likelihood problem is concave in linear utility coefficients \citep{McFadden1974,Train2009}. Its independence-of-irrelevant-alternatives (IIA) restriction, however, can be implausible when products naturally form groups. Adding a second economy fare should draw more demand from a similar fare than from a luxury product; adding a rail itinerary should substitute more strongly for another rail itinerary than for driving. The nested logit (NL) model captures this structure by partitioning alternatives into nests and assigning each nest a dissimilarity parameter that controls within-nest substitution \citep{Williams1977,DalyZachary1978,BenAkivaLerman1985}.

Nested logit has consequently become important in OR models of assortment and pricing. \citet{DavisGallegoTopaloglu2014} establish structural and approximation results for assortment optimization under several nested-logit variants. \citet{LiRusmevichientongTopaloglu2015} study assortment and price optimization for multilevel nesting, and \citet{chen2020capacitated} incorporate capacity constraints and endogenous prices. More recent exact algorithms solve broad nested-logit assortment variants by fractional programming and tailored branch-and-bound \citep{AlfandariHassanzadehLjubic2021}. This literature asks how to optimize an operational decision when the nested-logit parameters are known.

The corresponding estimation problem leaves a fundamental gap. Before a nested-logit demand model can support an operational decision, its utility coefficients and nest dissimilarities must be estimated jointly from observed choices. For fixed dissimilarities, the log-likelihood is concave in the utility coefficients under standard conditions \citep{DaganzoKusnic1993_convex_nested}. Joint estimation is nonconvex: dissimilarities divide utilities and multiply the inclusive values that aggregate within-nest attractiveness. Available implementations use Newton, quasi-Newton, sequential quadratic programming, or trust-region methods \citep{Bierlaire2003,Bierlaire2020,Brathwaite2018,NocedalWright2006}. Their stopping tests establish local stationarity or small steps, without bounding the distance in likelihood value from the global maximum. Multiple starts can improve an estimate but provide no global certificate. Consequently, even a globally optimized assortment or price vector may depend on an avoidably inferior local estimate of demand. The extensive optimization literature under nested-logit demand therefore calls for an accompanying global optimization foundation for estimating the model itself.

This gap raises two questions: can joint nested-logit estimation be solved with a verifiable global-accuracy guarantee, and which problem dimensions determine the cost of that guarantee? We answer both. We develop, to our knowledge, the first globally convergent framework with an explicit global performance guarantee for joint nested-logit maximum likelihood. The framework returns a feasible estimate together with upper and lower bounds on the globally optimal log-likelihood throughout the computation. On termination, its likelihood is within an additive tolerance $\varepsilon>0$ of the global maximum of the average log-likelihood over a declared compact parameter domain. We call this $\varepsilon$-global optimality; the guarantee concerns likelihood value rather than parameter distance.

Our central complexity result identifies the number of nests as the source of exponential computational cost. Let $N$ be the number of observations, $J$ the number of alternatives, $p$ the number of utility coefficients, and $L_f$ the number of free nest dissimilarity parameters. Let $C>0$ collect the scale-domain and uniform relaxation-error constants.
Under average normalization it is independent of $N$; any dependence on varying
model dimensions is assumed polynomial. Here $O(\cdot)$ suppresses constant factors and $\operatorname{poly}$ denotes a polynomial in its arguments. Under the compact-domain and quantitative relaxation-error conditions stated in Section~\ref{sec:global-bnb} and the appendix, the arithmetic-operation bound for additive average-likelihood error $\varepsilon$ has the form
\[
 O\!\left(
 \left(1+\frac{C}{\varepsilon}\right)^{L_f}
 \operatorname{poly}\!\left(N,J,p,\log\frac1\varepsilon\right)
 \right),
\]
Thus all exponential dependence is confined to the free nest count, while dependence on the observations, alternatives, and utility coefficients is polynomial. Fixing the number of nests therefore yields an additive polynomial-time approximation scheme (PTAS): for every fixed $\varepsilon>0$, an $\varepsilon$-globally optimal estimate is computable in time polynomial in the remaining input dimensions. This result is consequential for demand estimation because increasing the sample size or enriching the utility specification enlarges the polynomial part of the computation while preserving the dimension of the global search. The guarantee applies on the declared compact domain under the stated regularity and relaxation-error conditions.

The mechanism behind this result is a separation of convex estimation from nonconvex search. For fixed dissimilarities, we express maximum likelihood exactly as an exponential-cone program (ECP) \citep{Chares2009,MOSEKcookbook}; MNL follows by setting every dissimilarity to one. Our spatial branch-and-bound (B\&B) algorithm, called NL-B\&B, searches only over the free dissimilarity parameters \citep{HorstTuy1996,TawarmalaniSahinidis2002}. McCormick-envelope ECP relaxations supply valid likelihood upper bounds \citep{mccormick1976computability}, and fixed-dissimilarity ECPs supply feasible estimates and lower bounds. Every node requires polynomial work in the data and utility dimensions. The exponential factor arises from partitioning the free-dissimilarity domain, which explains both the fixed-nest PTAS and the practical value of a parsimonious nest structure.

Our main contribution is a global optimization foundation for jointly
estimating nested-logit utility coefficients and nest dissimilarities. The
novelty lies in a likelihood-specific conic decomposition and its accompanying
global-accuracy and complexity guarantees. Exponential-cone programming,
McCormick envelopes, and spatial B\&B provide the established tools used in
this construction. The contributions are:
\begin{itemize}
    \item \emph{Joint estimation with a global certificate.}
    We develop, to our knowledge, the first globally convergent framework with
    an explicit global-accuracy guarantee for joint nested-logit maximum
    likelihood. A feasible estimate and valid likelihood bounds quantify the
    remaining objective error throughout the search. Under the stated
    compact-domain and relaxation-error conditions, the method finitely
    attains any prescribed additive average-log-likelihood tolerance. This
    supplies an accuracy-based stopping criterion for an estimation problem
    commonly handled by local optimization (Section~\ref{sec:global-bnb}).

    \item \emph{An additive PTAS for a fixed number of nests.}
    We identify the number of free nest dissimilarities as the dimension
    governing exponential search complexity. For fixed nest count, the stated
    bounds give polynomial arithmetic complexity in the sample size,
    alternative count, utility dimension, and inverse target accuracy, under
    the declared domain and uniform relaxation-error assumptions. Thus the
    nonconvex estimation problem admits an additive PTAS in this regime.
    Increasing the amount of data enlarges the convex subproblems while
    preserving the dimension of the global search. A population-accuracy
    bound further connects the required sample size to the optimization
    tolerance (Section~\ref{sec:global-bnb}).

    \item \emph{An exact conic representation that enables the separation.}
    We represent fixed-dissimilarity nested-logit maximum likelihood exactly
    as an ECP and construct lifted conic relaxations for joint estimation.
    The fixed-scale programs yield feasible estimates, and the relaxations
    bound entire regions of the dissimilarity domain. Their combination
    allows the global method to branch only on free nest scales while
    optimizing all utility coefficients through convex programs. MNL follows
    as a special case, and a certified scale-grid method provides a direct
    alternative with an explicit discretization-error bound
    (Sections~\ref{sec:fixed-ecp}--\ref{sec:global-bnb}).

    \item \emph{Evidence that certification matters in practice.}
    On eight real discrete-choice datasets and 240 estimation instances,
    NL-B\&B certifies every reported estimate to the prescribed $1\%$
    relative tolerance. The comparisons expose inferior estimates returned
    by multistart local solvers despite normal termination: on the hardest
    common failure, even the best of five restarts of each of the three
    solvers falls $4.31\%$ below the certified reference. These results show
    both the feasibility of global certification on the benchmarks and the
    estimation errors that ordinary solver termination can leave undetected
    (Section~\ref{sec:results}).
\end{itemize}

Our contribution addresses a different stage of the analytics pipeline from
the literature on assortment, pricing, and revenue management under choice.
Those studies optimize an operational decision after a choice model and its
parameters have been specified. We address the upstream estimation problem:
given choice data and a nested-logit specification, we compute and certify a
globally near-optimal maximum-likelihood fit. Thus, existing results certify
the downstream decision conditional on the fitted demand model, whereas our
results certify the numerical quality of the parameter estimates supplied to
that decision model. The two guarantees are complementary and together reduce
the risk that a downstream policy is built on an avoidably poor local
likelihood solution.

The remainder of the paper is organized as follows.
Section~\ref{sec:literature} reviews the related choice-estimation and
operations literature. Section~\ref{sec:NL} introduces the nested-logit model
and its likelihood. Section~\ref{sec:fixed-ecp} develops the exact ECP for
fixed dissimilarity parameters and derives MNL as a corollary.
Section~\ref{sec:scale-grid} presents the direct certified grid method, and
Section~\ref{sec:global-bnb} develops the advanced spatial B\&B method together
with its convergence, complexity, and statistical guarantees.
Section~\ref{sec:results} reports the computational study, and
Section~\ref{sec:conclusion} concludes. The appendix contains the
complete proofs, B\&B formulation and implementation details, and additional
experiments.

\section{Related Literature}
\label{sec:literature}

\paragraph{Discrete-choice estimation.}
Random-utility foundations for MNL and generalized extreme-value models are
developed by \citet{McFadden1974}, \citet{Williams1977}, and
\citet{DalyZachary1978}; standard treatments include
\citet{BenAkivaLerman1985}, \citet{Train2009}, and \citet{HessDaly2014}.
Nested logit relaxes MNL's IIA restriction by allowing correlation among
alternatives in a prespecified hierarchy while retaining closed-form choice
probabilities. \citet{Daly1987} studies estimation of tree-logit models, and
\citet{DaganzoKusnic1993_convex_nested} establishes structural properties of
the fixed-scale likelihood. Broader members of the discrete-choice family
include mixed logit, which can approximate general random-utility models
\citep{McFaddenTrain2000}, and cross-nested and network-based specifications
that allow alternatives to share multiple substitution groups
\citep{VovshaBekhor2002,FosgerauMcFaddenBierlaire2013}. These richer models
increase flexibility but also make estimation and identification more
demanding. In practice, MNL and nested logit are commonly estimated with
Newton, quasi-Newton, or other local nonlinear routines implemented in
packages such as BIOGEME and PyLogit
\citep{Bierlaire2003,Bierlaire2020,Brathwaite2018}. Our focus is narrower but
fills a distinct gap: a global optimization certificate for the joint
maximum-likelihood estimation of the classical nested-logit model.

\paragraph{Assortment, pricing, and revenue management.}
Choice-based revenue management begins from the recognition that assortment,
availability, and price decisions alter substitution across products
\citep{Talluri2004revenue}. For MNL demand, tractable revenue-ordering and
related structural results support static, capacitated, and dynamic assortment
optimization \citep{rusmevichientong2010dynamic,rusmevichientong2014assortment}.
The nested-logit literature shows how much of this structure survives when
substitution is correlated within product groups. \citet{DavisGallegoTopaloglu2014}
analyze several nested-logit assortment variants, while
\citet{GallegoTopaloglu2014} treat constrained assortment and joint
assortment--pricing decisions. \citet{GallegoWang2014} derive pricing structure
and study multiproduct competition under product-differentiated price
sensitivities. \citet{LiRusmevichientongTopaloglu2015} extend assortment and
pricing results to multilevel nests, and \citet{chen2020capacitated} combine
capacity, assortment, and price decisions. Further work addresses display
positions \citep{ChenJiang2020} and exact algorithms under broad combinations
of nest coefficients and outside-option utilities
\citep{AlfandariHassanzadehLjubic2021}. All of these results take the demand
parameters as known; the quality of their operational prescriptions therefore
depends on the upstream estimation problem studied here.

\paragraph{Facility location and strategic-response models.}
Choice-based location models likewise embed probabilistic substitution into
long-lived siting decisions. Under nested-logit demand,
\citet{MendezVogelMarianovLuer2023} formulate the follower competitive facility
location problem, establish useful concavity, and develop an exact
branch-and-cut method. \citet{KimQiuShim2025} use nested logit to distinguish
private-car and park-and-ride alternatives in a capacitated hub-location
model. Moving to overlapping nests, \citet{LeEtAl2026CrossNestedLocation}
study competitive facility location under cross-nested logit and obtain an
exact outer-approximation approach through exponential-cone convexification.
Nested-logit structure also appears in strategic-response models. Logit
quantal response originated as a stochastic model of boundedly rational play
\citep{McKelveyPalfrey1995} and has since been embedded in Stackelberg security
games and algorithms for exploiting quantal opponents
\citep{yang2012computing,MilecEtAl2021}. Most directly,
\citet{MaiSinha2022NestedQR} introduce a nested-logit response rule for
correlated target choices in Stackelberg security games. Subsequent work uses
related quantal-response models for critical-service placement, network
interdiction, and heterogeneous adversary populations
\citep{MaiSinha2023Lifelines,MaiEtAl2024Interdiction,PhamEtAl2026MixedQR}.
These facility and game models reinforce the same point: errors in estimated
within-group substitution can propagate into locations, protection levels,
and equilibrium strategies.

\paragraph{Global and conic optimization.}
Spatial branch-and-bound constructs deterministic global bounds by recursively
tightening convex relaxations \citep{HorstTuy1996,TawarmalaniSahinidis2002}.
McCormick envelopes provide the convex hull of a bilinear graph over a box
\citep{mccormick1976computability}, while exponential cones represent
log-sum-exp constraints and support polynomial-time interior-point algorithms
\citep{Chares2009,BoydVandenberghe2004,NesterovNemirovskii1994}.
Conic representations have become especially useful when nonlinear choice
probabilities are embedded in mixed-integer optimization, as illustrated by
the cross-nested location model of \citet{LeEtAl2026CrossNestedLocation}. Our
formulation uses the same tools for a different object---the likelihood
itself. The global search is confined to nest dissimilarities, and every
utility-coefficient subproblem remains conic and convex.

\section{The Nested-Logit Model and Likelihood}
\label{sec:NL}
Nested logit groups correlated alternatives while retaining closed-form choice
probabilities \citep{Williams1977,DalyZachary1978,Train2009}. Its importance
comes from a central limitation of classical MNL: independence of
irrelevant alternatives forces the relative odds of two alternatives to remain
unchanged when similar alternatives are added or removed. Nested logit relaxes this
restriction parsimoniously by allowing alternatives in the same nest to share
unobserved factors and substitute more strongly with one another than with
alternatives in other nests.

This structure is directly relevant to operational demand models. Products may
share a brand, fare class, destination, technology, or service channel;
treating them as independent substitutes can distort predicted recapture and
therefore the optimal assortment, price vector, or capacity allocation
\citep{DavisGallegoTopaloglu2014,LiRusmevichientongTopaloglu2015,
chen2020capacitated}. Nested logit offers a useful middle ground: it is more
expressive than MNL, yet remains interpretable and far more parsimonious than
an unrestricted correlation model. Estimating the nest dissimilarities is
therefore operationally consequential because they determine the strength of
substitution within each product group.

We use $[K]=\{1,\ldots,K\}$ for any positive integer $K$ and
$|A|$ for the cardinality of a finite set $A$. There are $N$ observations,
$J$ alternatives in the universe $[J]$, and $p$ utility coefficients.
Observations are indexed by $n\in[N]$, alternatives by $j\in[J]$, nests by
$l\in[L]$, and coefficient coordinates by $d\in[p]$. Observation $n$ has a
nonempty available set (menu) $S_n\subseteq[J]$, an observed choice
$j_n\in S_n$, and a known attribute vector $\ba_{nj}\in\mathbb R^p$ for each
$j\in S_n$. The systematic utility is $v_{nj}=\bbt^\top\ba_{nj}$, where
$\bbt=(\beta_1,\ldots,\beta_p)^\top$ is the unknown coefficient vector and
$\top$ denotes transpose. The declared coefficient domain is a nonempty
compact set $\cB\subset\mathbb R^p$; the conic formulations below use a
convex, conically representable domain, such as the bounded coefficient box
used in the node relaxations.

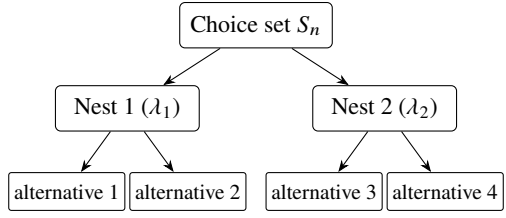
\begin{figure}[tb]
\centering
\begin{tikzpicture}[
  level/.style={draw,rounded corners=2pt,align=center,font=\footnotesize,
                minimum height=6mm,minimum width=19mm,inner xsep=4pt},
  alternative/.style={draw,rounded corners=1pt,align=center,font=\scriptsize,
                 minimum height=5mm,minimum width=12mm,inner xsep=2pt},
  >={Stealth[length=1.6mm]}]
\node[level] (root) at (0,0) {Choice set $S_n$};
\node[level] (n1) at (-1.7,-1.1) {Nest $1$ ($\lambda_1$)};
\node[level] (n2) at (1.7,-1.1) {Nest $2$ ($\lambda_2$)};
\node[alternative] (a1) at (-2.5,-2.2) {alternative $1$};
\node[alternative] (a2) at (-0.9,-2.2) {alternative $2$};
\node[alternative] (a3) at (0.9,-2.2) {alternative $3$};
\node[alternative] (a4) at (2.5,-2.2) {alternative $4$};
\draw[->] (root)--(n1);
\draw[->] (root)--(n2);
\draw[->] (n1)--(a1);
\draw[->] (n1)--(a2);
\draw[->] (n2)--(a3);
\draw[->] (n2)--(a4);
\end{tikzpicture}
\caption{Two-level nested-logit choice structure.}
\label{fig:nested-structure}
\end{figure}

Let the nonempty sets $\cN_1,\ldots,\cN_L$ partition $[J]$.
The active-nest set is $\cI_n=\{l\in[L]:\cN_l\cap S_n\neq\emptyset\}$,
and $l_n$ is the unique nest index with $j_n\in\cN_{l_n}$.
Nest $l$ has a dissimilarity (or scale) parameter
$\lambda_l\in[\underline\lambda_l,\overline\lambda_l]$, where
$0<\underline\lambda_l\leq\overline\lambda_l\leq1$ are declared bounds
and $\blambda=(\lambda_1,\ldots,\lambda_L)^\top$ is the full scale vector.
The usual unrestricted upper bound is $\overline\lambda_l=1$.
Smaller scales permit stronger within-nest dependence, whereas setting every
scale to one recovers MNL. Scales fixed in advance have equal lower and upper
bounds. In particular, we fix $\lambda_l=1$ for a singleton nest
($|\cN_l|=1$), since its scale cancels from the choice probabilities.
Let $\cF=\{l:\underline\lambda_l<\overline\lambda_l\}$ index the free
scales and let $L_f=|\cF|$. Only these coordinates are searched.
Figure~\ref{fig:nested-structure} illustrates the two-level choice structure.

\paragraph{Nested-logit likelihood.}
For an active nest $l\in\cI_n$, define the positive within-nest sum
\[
W_{nl}(\bbt,\lambda_l):=
\sum_{j\in\cN_l\cap S_n}
\exp\!\left(\frac{\bbt^\top\ba_{nj}}{\lambda_l}\right).
\]
This quantity aggregates the attractiveness of all available alternatives in
nest $l$. Throughout, $\log$ and $\ln$ denote the natural logarithm and $\exp$ its
inverse. The logarithm of $W_{nl}$ is the nest's inclusive value; the power
$W_{nl}^{\lambda_l}$ converts that inclusive value into attractiveness at the
upper level. The model can therefore be read as a two-stage choice. Write $P_n^{\mathrm{nest}}(l\mid S_n)$ for a nest-selection probability,
$P_n(j\mid l,S_n)$ for a within-nest probability, and $P_n(j\mid S_n)$ for
an alternative's unconditional choice probability. Their dependence on
$(\bbt,\blambda)$ is suppressed when clear. First, nest $l$ is selected with probability
\[
P_n^{\mathrm{nest}}(l\mid S_n)=
\frac{W_{nl}^{\lambda_l}}
{\sum_{r\in\cI_n}W_{nr}^{\lambda_r}}.
\]
Conditional on selecting $l$, alternative $j\in\cN_l\cap S_n$ is selected with
probability
\[
P_n(j\mid l,S_n)=
\frac{\exp(\bbt^\top\ba_{nj}/\lambda_l)}{W_{nl}}.
\]
Multiplying these two probabilities and taking logarithms gives the
contribution of the observed choice $j_n$ in nest $l_n$:
\[
\log P_n(j_n\mid S_n)=
(\lambda_{l_n}-1)\log W_{nl_n}
+\frac{\bbt^\top\ba_{nj_n}}{\lambda_{l_n}}
-\log\!\sum_{l\in\cI_n}W_{nl}^{\lambda_l}.
\]
The first term adjusts for within-nest substitution, the second is the scaled
utility of the chosen alternative, and the final term normalizes choice across
all active nests. Averaging these observation-level contributions yields the
joint maximum-likelihood estimation (MLE) problem, denoted NL-MLE:
\begin{equation}
\begin{aligned}
\max_{\bbt,\blambda}\quad
&\cL^{\NL}(\bbt,\blambda)\\
&=\frac1N\sum_{n=1}^N\Bigg[
 (\lambda_{l_n}-1)\log W_{nl_n}
 +\frac{\bbt^\top\ba_{nj_n}}{\lambda_{l_n}}
 -\log\!\sum_{l\in\cI_n}W_{nl}^{\lambda_l}\Bigg],\\
\mathrm{s.t.}\quad
&\bbt\in\cB,\quad
 \lambda_l\in[\underline\lambda_l,\overline\lambda_l],\quad l\in[L].
\end{aligned}
\label{prop:NL-MLE}\tag{\sf NL-MLE}
\end{equation}
Here $\cL^{\NL}$ is the average NL log-likelihood per observation.
All optimization objectives, primal--dual bounds, profile values, and additive
tolerances below use this average normalization. Multiplying the objective by
$N$ recovers the summed log-likelihood and leaves its maximizers unchanged. We write
$V^\star$ for its global maximum over the declared coefficient and scale
domains. For fixed $\blambda$, this objective is concave in $\bbt$
\citep{DaganzoKusnic1993_convex_nested}; joint estimation is nonconvex because
the scale parameters multiply inclusive values and divide utilities.  This
fixed-versus-joint distinction drives our algorithm: ECP solves each
fixed-scale problem globally, while branch-and-bound searches only over the
$L_f$ free scale parameters.

\section{ECP for Fixed Scale Parameters}
\label{sec:fixed-ecp}
This section isolates the tractable core of nested-logit estimation. We show
that once the dissimilarity vector is fixed, both levels of the likelihood can
be represented exactly with log-sum-exp epigraphs and exponential cones. The
resulting ECP is a globally solvable convex problem and later serves two roles:
it produces feasible estimates for the lower bound in NL-B\&B, and, when all
dissimilarities equal one, it immediately yields an ECP formulation and
complexity guarantee for MNL.

Fix an admissible $\blambda$. For each $n\in[N]$ and $l\in\cI_n$,
introduce a real epigraph variable $z_{nl}\geq\log W_{nl}$ and a real variable
$y_n\geq\log\sum_{l\in\cI_n}e^{\lambda_lz_{nl}}$. An epigraph variable
bounds a function from above. These inequalities can be made tight at an
optimum without decreasing the likelihood. Introduce nonnegative variables
$k_{njl}$ for $j\in\cN_l\cap S_n$ and $h_{nl}$ for $l\in\cI_n$ to bound,
respectively, $e^{v_{nj}/\lambda_l-z_{nl}}$ and
$e^{\lambda_lz_{nl}-y_n}$. Bold symbols $\bz,\by,\bk,\bh$ collect all
corresponding components; variables with nest indices exist only for active
observation--nest pairs. Thus $\bz,\bh\in\mathbb R^\Lambda$,
$\by\in\mathbb R^N$, and $\bk\in\mathbb R^Z$, where
$\Lambda=\sum_{n=1}^N|\cI_n|$ and $Z=\sum_{n=1}^N|S_n|$ count active
observation--nest and observation--alternative pairs, respectively.
Denote this fixed-scale formulation by NL-ECP.

We use the exponential cone
\[
\cK_{\exp}:=\operatorname{cl}\{(a,b,c)\in\mathbb R^3:
 b>0,\ a\geq b\exp(c/b)\},
\]
where $\operatorname{cl}$ denotes closure. This coordinate order places the
exponential epigraph value first, so $(a,1,c)\in\cK_{\exp}$ means
$a\geq e^c$ \citep{Chares2009,MOSEKcookbook}.

\begin{proposition}[Fixed-scale nested-logit ECP]
\label{prop:fixed-nl-ecp}
For any fixed admissible $\blambda\in(0,1]^L$, nested-logit maximum likelihood is equivalent
to
\begin{align}
   \max_{\bbt\in\cB, \bz, \by, \bk, \bh} \quad
    & \frac1N\sum_{n=1}^N \Big[ (\lambda_{l_n}-1) z_{n l_n} \label{prob:NL-ECP}\tag{\sf NL-ECP}\\
    & \qquad\quad + \tfrac{\bbt^\top \ba_{n j_n}}{\lambda_{l_n}} - y_n \Big], \nonumber\\
    \text{s.t.}\quad & \sum_{j\in \cN_l\cap S_n} k_{njl} \leq 1,
    \quad \forall n,\ l\in\cI_n, \nonumber\\
    & \sum_{l\in \cI_n} h_{nl} \leq 1,
    \quad \forall n,\nonumber \\
    & (k_{njl},1,\bbt^\top\ba_{nj}/\lambda_l-z_{nl})\in\cK_{\exp},
    \nonumber\\[-2pt]
    & \hspace{22mm}\forall n,\ l\in\cI_n,\ j\in\cN_l\cap S_n, \nonumber\\
    & (h_{nl},1,\lambda_lz_{nl}-y_n)\in\cK_{\exp},
    \quad \forall n,\ l\in\cI_n.\nonumber
\end{align}
\end{proposition}

\paragraph{Proof idea and implications.}
The proof introduces one epigraph variable for each within-nest log-sum-exp and
one for the upper-level log-sum-exp. Each exponential term is represented by a
three-dimensional exponential-cone constraint. Because the epigraph variables
enter the maximization objective with the appropriate nonpositive
coefficients, the inequalities are tight at an optimum, so the conic program
neither relaxes nor restricts the fixed-$\blambda$ likelihood. The result is the
key structural decomposition behind the paper: fixing only the nest
dissimilarities turns all remaining estimation variables into a convex,
globally solvable subproblem with a primal--dual certificate. The full proof appears in the appendix, p.~\pageref{ec:proof:fixed-nl-ecp}.

\begin{corollary}[MNL as a special case]
\label{cor:mnl-ecp}
Setting $\lambda_l=1$ for every nest collapses nested logit to MNL. After
eliminating the redundant nest-level variables from \eqref{prob:NL-ECP}, the
MNL maximum-likelihood problem, denoted MNL-ECP, is exactly the following.
Here $t_n\in\mathbb R$ bounds the observation-level log-sum-exp from above,
$u_{nj}\geq0$ bounds $e^{v_{nj}-t_n}$, and $\bt,\bu$ collect these variables:
\begin{align}
 \max_{\bbt\in\cB,\bt,\bu}\quad
 &\frac1N\sum_{n=1}^N\big(\bbt^\top\ba_{nj_n}-t_n\big)
 \label{prob:MNL-ECP}\tag{\sf MNL-ECP}\\
 \mathrm{s.t.}\quad
 &\sum_{j\in S_n}u_{nj}\leq1, && n\in[N],\nonumber\\
 &(u_{nj},1,\bbt^\top\ba_{nj}-t_n)\in\cK_{\exp},
 &&n\in[N],\ j\in S_n.\nonumber
\end{align}
Recall the count $Z$ of available observation--alternative pairs, and let $\varepsilon>0$ be the desired
additive objective accuracy. A path-following primal--dual interior-point method
\citep{NesterovNemirovskii1994,BenTalNemirovski2001}
computes an $\varepsilon$-optimal MNL estimate in
\[
 O\!\left(\sqrt{Z}\log(1/\varepsilon)
 \left[Zp^2+(p+N)^3\right]\right)
\]
arithmetic operations. Thus MNL estimation is polynomial in the sample,
alternative, and feature dimensions and logarithmic in the requested conic
accuracy.
\end{corollary}

\paragraph{Proof idea and implications.}
Setting every $\lambda_l$ to one makes the within-nest and across-nest
log-sum-exp expressions compose into a single log-sum-exp over $S_n$; removing
the redundant nest variables gives MNL-ECP directly. The complexity bound
combines the $O(\sqrt Z\log(1/\varepsilon))$ iteration bound for a product of
$Z$ exponential cones with a conservative bound on the Newton linear algebra.
This shows that MNL requires no separate global-search machinery and clarifies
the source of difficulty in nested logit: nonconvexity arises from estimating
the dissimilarities, not from the utility coefficients or the sample size. The full proof appears in the appendix, p.~\pageref{ec:proof:mnl-ecp}.

\paragraph{Advantages of the MNL-ECP reformulation.}
A standard Newton method works directly with the $p$ utility coefficients. One
iteration evaluates all $Z$ alternative occurrences, forms a Hessian in
$O(Zp^2)$ operations, and solves a $p\times p$ system in $O(p^3)$ dense
arithmetic. It can therefore be substantially faster than MNL-ECP, which adds
$N+Z$ auxiliary variables and whose displayed worst-case bound contains the
conservative term $(p+N)^3$. The dependence on the target accuracy requires
care. Without uniform curvature, a Lipschitz Hessian, and a specified
globalization rule, an ordinary Newton implementation has no uniform
$\varepsilon$-complexity guarantee. Under the standard assumptions that the
negative log-likelihood is uniformly strongly convex on the parameter domain,
its Hessian is Lipschitz, and damped Newton steps are chosen by a suitable line
search, a conservative global linear-convergence bound gives
\[
 O\!\left(\left[Zp^2+p^3\right]\log\frac1\varepsilon\right)
\]
arithmetic operations to reach objective error $\varepsilon$. Once the iterates
enter the quadratic-convergence neighborhood, the remaining accuracy requires
only $O(\log\log(1/\varepsilon))$ Newton steps, or
$O([Zp^2+p^3]\log\log(1/\varepsilon))$ additional arithmetic operations
\citep{NocedalWright2006}. These rates hide constants depending on conditioning
and on the size of the region needed to reach the local quadratic phase.

The advantage of the ECP
is not a claim of uniformly lower runtime. Rather, the formulation can be sent
unchanged to mature conic solvers, requires no hand-derived gradient, Hessian,
line search, or regularization strategy, and returns primal and dual residuals
together with a solver-checkable duality gap \citep{MOSEKcookbook}. Its
$O(\sqrt Z\log(1/\varepsilon))$ path-following iteration bound also gives an
explicit polynomial worst-case guarantee, whereas the practical convergence
of an ordinary Newton implementation depends on conditioning, globalization,
and its stopping rule. Most importantly here, MNL-ECP uses the same conic
language as fixed-scale NL-ECP and the node relaxations, so one solver
infrastructure supports the MNL benchmark, feasible nested-logit incumbents,
and global upper bounds.

\section{Direct Global Estimation by Scale Discretization}
\label{sec:scale-grid}

A direct way to exploit Proposition~\ref{prop:fixed-nl-ecp} is to discretize
the low-dimensional domain of the free dissimilarity parameters and solve one
NL-ECP at every grid point. This method requires no nonconvex solver: all
utility coefficients are optimized globally within each fixed-scale problem,
and the best grid solution is selected. Besides providing a transparent
baseline for global estimation, the construction exposes why the number of
nests, rather than the number of observations or alternatives, is the critical
nonconvex dimension.

Recall the free-scale set $\cF$ and its size $L_f$. Define the free-scale box
$\cD=\prod_{l\in\cF}[\underline\lambda_l,\overline\lambda_l]$.
Write $\blambda_{\cF}$ for the vector of free coordinates. When writing
$\blambda\in\cD$, we mean that $\blambda_{\cF}\in\cD$ and that all other
coordinates retain their declared fixed values. This convention also applies
to grids in $\cD$ and to maximizations over this box. Define the interval
length $\Delta_l=\overline\lambda_l-\underline\lambda_l>0$ for $l\in\cF$.
We call the resulting direct grid estimator NL-GRID.
The profile log-likelihood is
\begin{equation}
 V(\blambda):=\max_{\bbt\in\cB}\cL^{\NL}(\bbt,\blambda)
 =\operatorname{val}\big(\text{NL-ECP}(\blambda)\big).
 \label{eq:profile-likelihood}
\end{equation}
Here $\operatorname{val}$ denotes the optimal objective value of a program.
For a target mesh width $h>0$, set
$K_l=\lceil(\overline\lambda_l-\underline\lambda_l)/h\rceil$ (where $\lceil x\rceil$ is the least integer not smaller than $x$) and form the
endpoint-inclusive grid
\begin{equation}
 \cG_h=\mathop{\times}_{l\in\cF}
 \left\{\underline\lambda_l+
 \frac{i}{K_l}(\overline\lambda_l-\underline\lambda_l):
 i=0,\ldots,K_l\right\}.
 \label{eq:scale-grid}
\end{equation}
For a conic objective tolerance $\eta>0$, at every $\blambda\in\cG_h$, solve \eqref{prob:NL-ECP} to conic duality gap
at most $\eta$. Let $\underline V_\eta(\blambda)$ and $\overline V_\eta(\blambda)$ be its
primal lower and dual upper objective bounds. Using $\arg\max$ for the set of
maximizers, define the direct grid estimator and its lower bound $LB_h$ by
\begin{equation}
 \widehat\blambda_h
 \in\argmax_{\blambda\in\cG_h}
 \underline V_\eta(\blambda),
 \qquad
 LB_h:=\max_{\blambda\in\cG_h}\underline V_\eta(\blambda).
 \label{prob:NL-GRID}\tag{\sf NL-GRID}
\end{equation}
Let $\widehat\bbt_h$ be the feasible coefficient vector returned by the
primal solve at $\widehat\blambda_h$. Thus the maximization selects a scale
vector, while its associated conic solution supplies the coefficient estimate.
The upper bound $UB_h$ is defined in the following proposition.
All grid-point ECPs are independent and can therefore be solved in parallel.

\begin{proposition}[Grid error bound]
\label{prop:grid-error}
Suppose $\cB\times\cD$ is compact and $\underline\lambda_l>0$ for every free
scale. Let $h$ be the maximum coordinate-wise grid spacing, let $\eta$ be the
maximum primal--dual gap of each grid-point NL-ECP, and let $G_\lambda$ be any
uniform scale-sensitivity bound satisfying (with $\nabla_{\blambda_{\cF}}$
denoting the gradient in the free scales and $\|x\|_1=\sum_i|x_i|$)
\[
 G_\lambda\ \geq\
 \sup_{(\bbt,\blambda)\in\cB\times\cD}
 \left\|\nabla_{\blambda_{\cF}}
 \cL^{\NL}(\bbt,\blambda)\right\|_1.
\]
Then the output of \eqref{prob:NL-GRID} satisfies
\[
 0\leq V^\star-LB_h\leq \frac{G_\lambda h}{2}+\eta,
 \qquad
 V^\star:=\max_{\blambda\in\cD}V(\blambda).
\]
Moreover,
\[
 UB_h:=\max_{\blambda\in\cG_h}
 \overline V_\eta(\blambda)+\frac{G_\lambda h}{2}
\]
is a valid global upper bound, so $[LB_h,UB_h]$ is an explicit optimality
certificate. If $G_\lambda>0$, choosing
$h=\varepsilon/G_\lambda$ and $\eta=\varepsilon/2$ yields an
$\varepsilon$-globally optimal nested-logit estimate.
\end{proposition}

The full proof appears in the appendix, p.~\pageref{ec:proof:grid-error}.

\paragraph{Proof idea and implications.}
The derivative bound and the mean-value theorem imply that
$\cL^{\NL}(\bbt,\cdot)$ is $G_\lambda$-Lipschitz in the infinity norm $\|x\|_\infty=\max_i|x_i|$,
uniformly over $\bbt\in\cB$. Maximization over the common set $\cB$ preserves
this Lipschitz constant, so the profile likelihood $V$ is also
$G_\lambda$-Lipschitz. Every point of $\cD$ lies within infinity-norm distance
$h/2$ of an endpoint-inclusive grid point. Applying this fact to a global
maximizer and adding the conic duality gap proves both bounds. Thus the grid
search is a genuinely certified approximation method when a valid
$G_\lambda$ is available; such a constant can be bounded before optimization
from the declared coefficient, attribute, and positive scale bounds. If
$G_\lambda=0$, the profile likelihood is constant over $\cD$, so one scale
point suffices.

\paragraph{Overall complexity.}
Recall the counts $Z$ and $\Lambda$ and the free-scale interval lengths
$\Delta_l$.

\begin{proposition}[Complexity of scale-grid estimation]
\label{prop:grid-complexity}
For mesh width $h>0$ and conic accuracy $\eta>0$, the sequential arithmetic
complexity of NL-GRID is
\begin{equation}
 O\!\left(
 \prod_{l\in\cF}
 \left(1+\left\lceil\frac{\Delta_l}{h}\right\rceil\right)
 \sqrt{Z+\Lambda}\log\frac1\eta
 \left[Zp^2+(p+N+\Lambda)^3\right]
 \right).
 \label{eq:grid-complexity-general}
\end{equation}
In particular, the choices $h=\varepsilon/G_\lambda$ and
$\eta=\varepsilon/2$ from Proposition~\ref{prop:grid-error} give
\begin{equation}
 O\!\left(
 \prod_{l\in\cF}
 \left(1+\left\lceil\frac{G_\lambda\Delta_l}{\varepsilon}\right\rceil\right)
 \sqrt{Z+\Lambda}\log\frac1\varepsilon
 \left[Zp^2+(p+N+\Lambda)^3\right]
 \right).
 \label{eq:grid-complexity}
\end{equation}
Thus NL-GRID is polynomial in the data dimensions and $1/\varepsilon$ for
fixed $L_f$, while its dependence on the number of free nests is exponential.
\end{proposition}

\paragraph{Proof idea and implications.}
Coordinate $l$ contributes
$1+\lceil\Delta_l/h\rceil$ grid values, so their Cartesian product gives the
leading factor in \eqref{eq:grid-complexity-general}. Each point requires one
fixed-scale NL-ECP with $Z+\Lambda$ exponential-cone blocks; multiplying its
interior-point complexity by the grid cardinality gives the result. The ECPs
can be solved in parallel, but this reduces wall-clock time rather than total
work. The bound makes the limitation of direct discretization explicit: it
solves every grid point, including regions that are visibly suboptimal.
Section~\ref{sec:global-bnb} instead uses valid regional upper bounds to prune
such regions while retaining the same global guarantee. The full proof appears in the appendix, p.~\pageref{ec:proof:grid-complexity}.

\section{Advanced Spatial B\&B for Global Nested-Logit Estimation}
\label{sec:global-bnb}

The central idea for obtaining a globally optimal nested-logit estimate---and,
in particular, an additive PTAS when the number of nests is fixed---is to
avoid the exhaustive grid search of Section~\ref{sec:scale-grid} by bounding
entire boxes of scale vectors at once. We separate the small nonconvex part of
the estimation problem from its large convex part. The dissimilarity parameters are the
only variables responsible for nonconvexity; once they are fixed, the utility
coefficients are obtained globally from NL-ECP. We therefore partition only
the dissimilarity domain. For each partition box, a convex conic relaxation
optimistically bounds the best likelihood attainable anywhere in that box,
while fixing a candidate dissimilarity vector and solving NL-ECP produces an
attainable likelihood. The former is a global upper bound and the latter a
global lower bound. Boxes that cannot improve the best attainable solution are
discarded, and the remaining boxes are successively refined until the upper
and lower bounds differ by at most the prescribed tolerance. This mechanism
provides both an anytime certificate during the search and an
$\varepsilon$-global-optimality certificate at termination. This follows the
standard spatial branch-and-bound paradigm of bounding, pruning, and exhaustive
refinement \citep{HorstTuy1996,TawarmalaniSahinidis2002}.

Joint estimation is nonconvex because NL-ECP contains the products
$\lambda_lz_{nl}$ and the ratios $\bbt^\top\ba_{nj}/\lambda_l$.  Our spatial
branch-and-bound algorithm NL-B\&B searches only over the $L_f$ free scale parameters.
At a node indexed by $r$, let $\underline\lambda_l^r$ and
$\overline\lambda_l^r$ denote its lower and upper scale bounds, so
$\lambda_l\in[\underline\lambda_l^r,\overline\lambda_l^r]$.  We lift
$q_{nl}=\lambda_lz_{nl}$, $\rho_l=1/\lambda_l$, and
$\bgamma_l=\rho_l\bbt$, then replace the bilinear equalities by McCormick
envelopes over the current box \citep{mccormick1976computability}. Together with the remaining exponential-cone
constraints, these envelopes form a convex ECP relaxation whose value is a
valid node upper bound. Here $q_{nl}\in\mathbb R$ is a lifted product,
$\rho_l>0$ a reciprocal scale, and $\bgamma_l\in\mathbb R^p$ a vector of
scaled coefficients. Node superscripts specify interval bounds, not powers;
$\underline\lambda_l^r,\overline\lambda_l^r$ are the node's scale endpoints.
Let $UB_r$ denote a valid node upper bound, $LB$ the best feasible likelihood
found so far, and $UB$ the global upper bound maintained across nodes.
The global absolute gap is $UB-LB$. With $\cA$ denoting the open-node queue,
a conservative certificate retaining tolerance-pruned regions is
\[
 UB=\max\{LB+\varepsilon,\ \max_{r\in\cA}UB_r\},
 \qquad LB\leq V^\star\leq UB,
\]
where the maximum over an empty queue is $-\infty$. The tolerance
$\varepsilon>0$ here is additive in average log-likelihood per observation. Hats indicate
computed estimates or candidate values; a superscript $\mathrm{inc}$ identifies
the current best feasible estimate (the incumbent), and a star identifies an
exact optimizer or optimal value.

For a refinement threshold $\delta>0$, a free scale is eligible for branching
if its node interval has width greater than $\delta$ or endpoint ratio greater
than $1+\delta$.
A lower bound is obtained by fixing $\blambda$ to a candidate from the
relaxation and solving the exact fixed-scale NL-ECP for $\bbt$.  Nodes whose
upper bounds cannot improve the incumbent are pruned; otherwise the algorithm
branches on an eligible scale interval contributing the largest relaxation error.  As
the scale boxes shrink, the McCormick envelopes converge to the lifted
equalities, so the upper and lower bounds meet.  The method is therefore
anytime: it returns a feasible estimate and a valid global optimality gap at
interruption, and certifies the requested tolerance upon termination.  The full
node relaxation, tightening and branching rules, and convergence argument
appear in the appendix under ``Full B\&B Details.''

For the quantitative convergence and complexity bounds, assume the following
uniform relaxation-error condition. There is a constant $H>0$, independent of
$N$, such that a node whose free-scale intervals all have width at most
$\delta>0$ and endpoint ratio at most $1+\delta$ has a certified upper bound
within $H\delta+\eta$ of its fixed-scale feasible probe. Here $\eta>0$ is
the combined error allowance for the relaxation and probe solves; each solve
receives at most $\eta/2$. The condition is imposed on the declared bounded
domains, and dependence of $H$ on other varying problem dimensions must be
polynomial for the corresponding polynomial-time claim. Choose
$H\delta+\eta\leq\varepsilon$. The global stopping test is always
$UB-LB\leq\varepsilon$; the interval conditions are a sufficient refinement
threshold for this test, not an alternative certificate.

Algorithm~\ref{alg:global-NL} summarizes the certified B\&B procedure used in
our analysis and experiments. Its two bounding steps use the same NL-ECP
machinery: the McCormick--ECP relaxation supplies a valid upper bound, whereas
fixing the dissimilarity parameters and solving NL-ECP supplies a feasible
lower bound. Under the relaxation-error condition, a box with no eligible
scale satisfies $UB_r-LB\leq\varepsilon$ after its probe and is therefore
pruned by the same objective-gap test.

\begin{algorithm}[tb]
\caption{B\&B for Globally Certified Nested-Logit Estimation}
\label{alg:global-NL}
\begin{algorithmic}[1]
\REQUIRE Domains $\cB,\cD$, average-objective tolerance $\varepsilon>0$,
refinement threshold $\delta>0$, and conic allowance $\eta>0$ with
$H\delta+\eta\leq\varepsilon$; a time limit.
\STATE Obtain a feasible incumbent
$(\widehat\bbt^{\mathrm{inc}},\widehat\blambda^{\mathrm{inc}})$ and its
average log-likelihood $LB$ from fixed-scale NL-ECP solves.
\STATE Initialize the best-first queue $\cA$ with the root box and a valid
upper bound; set $UB=\max\{LB+\varepsilon,\max_{r\in\cA}UB_r\}$.
\WHILE{$UB-LB>\varepsilon$, the queue is nonempty, and time remains}
  \STATE Pop the highest-bound box $r$ and tighten its bounds. If
  $UB_r>LB+\varepsilon$, solve its relaxation for a certified upper bound and
  fix a candidate scale to obtain a feasible probe; use conic allowances
  $\eta/2$ for each solve and update $LB$ and the incumbent.
  \STATE If $UB_r\leq LB+\varepsilon$, prune $r$. Otherwise branch on a free
  scale whose width exceeds $\delta$ or whose endpoint ratio exceeds
  $1+\delta$, choosing the largest score among these eligible scales;
  insert viable children with valid upper bounds.
  \STATE If time expires during a node operation, reinsert the unresolved
  node with its last valid bound. Recompute
  $UB=\max\{LB+\varepsilon,\max_{r\in\cA}UB_r\}$.
\ENDWHILE
\RETURN The incumbent and $LB\leq V^\star\leq UB$; report whether
$UB-LB\leq\varepsilon$ was attained or the run stopped at its time limit.
\end{algorithmic}
\end{algorithm}

{
\begin{theorem}[Validity and convergence of NL-B\&B]
\label{thm:bnb-validity}
Suppose $\cB$ is compact, every scale is bounded away from zero, and all
McCormick bounds, tightening steps, and cuts are valid. Use certified dual
upper bounds and feasible-primal lower bounds, with the conic error allowances
included conservatively in the reported bounds. Then at every interruption,
$LB\leq V^\star\leq UB$ for the bounded average-likelihood problem.
Under the uniform relaxation-error condition above, choose
$H\delta+\eta\leq\varepsilon$ and branch only on intervals that have not met
both refinement conditions. In the absence of a time limit,
Algorithm~\ref{alg:global-NL} terminates after finitely many nodes when
$UB-LB\leq\varepsilon$ and returns an $\varepsilon$-globally optimal estimate
in average log-likelihood. A time-limited run retains its certificate without
claiming the requested tolerance unless this gap test is satisfied.
\end{theorem}}

\paragraph{Proof idea and implications.}
At every node, McCormick envelopes contain the lifted bilinear graph, so the
node ECP is a valid relaxation and its value is an upper bound. Conversely,
fixing a candidate dissimilarity vector and solving NL-ECP produces a feasible
estimate and hence a lower bound. Branching refines only eligible intervals. The uniform error condition
ensures that a box meeting both interval thresholds is prunable by the
objective-gap test. The finite refinement count for the declared positive
scale bounds therefore gives a finite search for any $\varepsilon>0$. The
theorem therefore upgrades branch-and-bound from a heuristic search into an
anytime estimator: even an interrupted run reports a valid global gap, and a
completed run certifies the requested accuracy.
The full proof appears in the appendix, p.~\pageref{ec:proof:bnb-validity}.

%\textcolor{black}{The theorem concerns the explicitly bounded estimationproblem. It does not claim that an arbitrary finite box contains an unconstrained MLE; the appendix explains how the bounds are selected and stress-tested by systematic expansion.}

\subsection{Complexity and Practical Implications}
Let $Z=\sum_n|S_n|$ be the number of active observation--alternative pairs and
$\Lambda=\sum_n|\cI_n|$ the number of active observation--nest pairs.  A node
contains $Z+\Lambda$ exponential-cone blocks and polynomially many lifted
variables and linear inequalities.  Hence all work at a fixed node---including
the relaxation and fixed-scale incumbent solve---is polynomial in the data and
model dimensions.  The non-polynomial factor comes only from partitioning the
$L_f$-dimensional free-scale box.

\begin{proposition}[Nested-logit complexity]
\label{prop:nl-bnb-complexity-main}
Let $\delta>0$ be the scale-box refinement tolerance, $\eta>0$ the conic-solver
combined accuracy allowance, and assume the uniform relaxation-error
condition with $H\delta+\eta\leq\varepsilon$. Let $D$ be the dimension of the global Schur-complement
linear system in a conic Newton step. Define
\begin{align*}
 \kappa_l&=\max\!\left\{\overline\lambda_l-\underline\lambda_l,
 \log\!\frac{\overline\lambda_l}{\underline\lambda_l}\right\},
 &\kappa&=\max_{l\in[L]}\kappa_l,\\
 D&=O\!\left(p+Lp+L+N+\Lambda\right).
\end{align*}
Let $C_{\mathrm{node}}$ and $C_{\mathrm{total}}$ denote arithmetic-operation
counts for one node and for the full search, respectively. These costs are,
respectively,
\begin{align*}
 C_{\mathrm{node}}
   &=O\!\left(\sqrt{Z+\Lambda}\,\log\frac{1}{\eta}
      \left[Zp^2+D^3\right]\right),\\
 C_{\mathrm{total}}
   &=O\!\left(\left(1+\frac{\kappa}{\delta}\right)^{{L_f}}
      C_{\mathrm{node}}\right).
\end{align*}
{Here the exponent is the number $L_f$ of free scales; fixed and
singleton-nest scales do not create branching coordinates. Thus, for fixed
$L_f$, a fixed compact domain, and fixed tolerances,} the algorithm is polynomial in the
numbers of observations, alternatives, and utility parameters; its exponential
dependence is confined to the number of nests.
\end{proposition}

\paragraph{Proof idea and implications.}
The proof first counts $Z+\Lambda$ exponential-cone blocks and the lifted
variables in one node, then applies the standard interior-point iteration and
linear-algebra bounds \citep{NesterovNemirovskii1994,BenTalNemirovski2001}. It separately counts how many arithmetic or geometric
subintervals of width $\delta$ can cover each dissimilarity range; taking the
Cartesian product over the $L_f$ free scales yields
$(1+\kappa/\delta)^{L_f}$. Thus the bound separates polynomial per-node ECP cost from an exponential search
over the $L_f$ free scales. Increasing the sample size enlarges each convex node
problem but does not add a branching coordinate. This fixed-parameter structure
is favorable in applications, where the number of nests is typically small;
the appendix gives the complete counting argument. In particular, choose
$\delta$ and $\eta$ so that the relaxation and conic errors sum to at most a
requested additive average-likelihood tolerance $\varepsilon$. For every fixed $L_f$
and fixed compact parameter domain, Proposition~\ref{prop:nl-bnb-complexity-main}
then gives a polynomial-time approximation scheme: for each fixed
$\varepsilon>0$, the running time is polynomial in $N$, $Z$, and $p$, and the
returned likelihood is within $\varepsilon$ of the global maximum. We call the
scheme additive because log-likelihoods can be negative, making a multiplicative
ratio inappropriate. The constants deteriorate exponentially with $L_f$, as is
explicit in the bound. The full proof can be found in the appendix, p.~\pageref{ec:proof:nl-bnb-complexity}.

\subsection{Sample Efficiency and Computational Cost}
A globally optimal empirical estimate may still generalize poorly when learned
from too few choices. This is especially relevant for nested logit, whose scale
parameters capture correlation but increase the quantities to be learned.
Larger samples reduce statistical error while enlarging each likelihood ECP.
The following result connects the observations needed for population accuracy
to the cost of solving the empirical problem.

\begin{proposition}[End-to-end accuracy]
\label{prop:sample-computation}
Let $\theta=(\bbt,\blambda_{\cF})\in\Theta\subset\mathbb R^{p+L_f}$
collect the unknown coefficients and free scales. Let $X$ be a random choice
observation (menu, attributes, and selected alternative), and let
$\ell(\theta;X)$ be its NL log-likelihood. For independent and identically
distributed (i.i.d.) copies $X_1,\ldots,X_N$, define the population and empirical average log-likelihoods
$Q(\theta)=\mathrm E[\ell(\theta;X)]$ and
$\widehat Q_N(\theta)=N^{-1}\sum_{n=1}^N\ell(\theta;X_n)$. Here
$\mathrm E$ denotes expectation and $\Pr$ denotes probability under the
observation law. Under the conditions in the appendix, let
$\Theta$ be compact with Euclidean diameter
$R=\sup_{\theta,\theta'\in\Theta}\|\theta-\theta'\|_2$, where
$\|x\|_2=(\sum_i x_i^2)^{1/2}$. Let $G>0$ and $B>0$ satisfy, for all $\theta,\theta'\in\Theta$,
$|\ell(\theta;x)-\ell(\theta';x)|\leq G\|\theta-\theta'\|_2$ for every
possible observation $x$, and let the range of $\ell(\theta;X)$ have length
at most $B$ for every $\theta$. Let $\theta_0\in\arg\max_{\theta\in\Theta}Q(\theta)$
be a population maximizer, and $\widehat\theta_N\in\Theta$ the computed
empirical estimate. For target population error $\varepsilon_{\mathrm{stat}}>0$
and failure probability $\alpha\in(0,1)$, if, for a universal constant $c>0$,
\[
 N\geq c\frac{B^2}{\varepsilon_{\mathrm{stat}}^2}
 \left[(p+L_f)\log\!\left(1+\frac{GR}{\varepsilon_{\mathrm{stat}}}\right)
 +\log\frac1\alpha\right]
\]
and the empirical problem is solved to average log-likelihood error at most
$\varepsilon_{\mathrm{stat}}/3$, then
\[
 \Pr\!\left\{Q(\theta_0)-
 Q(\widehat\theta_N)\leq\varepsilon_{\mathrm{stat}}\right\}\geq1-\alpha.
\]
Thus the statistical and optimization errors can be matched at the same target
accuracy.
\end{proposition}

\paragraph{Proof idea and implications.}
The proof uses a finite $\varepsilon_{\mathrm{stat}}$-net of the compact parameter set and a
uniform concentration bound for the empirical log-likelihood on that net.
Lipschitz continuity extends the bound from the net to the full parameter
space, and the standard empirical-risk decomposition adds two uniform sampling
errors and the certified optimization error. The implication is operational:
global optimization alone does not remove sampling error, so the solver
tolerance should be chosen on the same scale as statistical precision rather
than made arbitrarily small. The full proof appears in the appendix, p.~\pageref{ec:proof:sample-computation}.

{
\begin{corollary}[Computational consequence]
\label{cor:sample-computation}
Under the conditions of Proposition~\ref{prop:sample-computation} and bounded
menus, suppose the NL node relaxation satisfies the quantitative error
condition stated in the appendix. For fixed $L_f$ and a fixed compact
parameter domain, NL-B\&B attains the statistical guarantee in arithmetic time
polynomial in $N$, $p$, and $1/\varepsilon_{\mathrm{stat}}$; its exponential dependence is
confined to $L_f$.
\end{corollary}}

\paragraph{Proof idea and implications.}
The corollary substitutes the sample-size requirement of
Proposition~\ref{prop:sample-computation} and an $O(\varepsilon_{\mathrm{stat}})$ empirical
optimization tolerance into the NL-B\&B arithmetic bound. For
fixed $L_f$, the resulting cost is polynomial in the sample and feature
dimensions and in $1/\varepsilon_{\mathrm{stat}}$; only the nest count remains in the exponent.
Consequently, collecting more observations makes each conic problem larger but
does not make the global-search dimension larger, which is the central reason
the method can scale statistically when the nesting structure is small. The full proof can be found in the appendix, p.~\pageref{ec:proof:sample-computation-cor}.

%A proof, the shared quantitative relaxation-error assumption, and an explicitnested-logit arithmetic bound appear in the appendix section ``Sample Efficiency and Computational Complexity.''

Proposition~\ref{prop:sample-computation} controls excess expected
log-likelihood, the population counterpart of the estimation objective; it does
not by itself assert parameter-distance or choice-probability error. The result
has two practical implications. First, collecting more data
increases the size of the convex node problems but not the dimension of the
combinatorial search, so large samples remain compatible with global
certification when the nest structure is small. Second, the certified solver
tolerance should be matched to statistical precision: solving far below the
$O(N^{-1/2})$ sampling error adds computation with little population-level
benefit. Thus one can choose both $N$ and the optimality tolerance from a target
population log-likelihood accuracy, while using the smallest behaviorally supported number of
nests. The appendix section ``Sample Efficiency and Computational Complexity''
gives the assumptions and the nested-logit end-to-end bound.

\paragraph{Comparison with the scale-grid method.}
In theory, NL-GRID in Section~\ref{sec:scale-grid} and NL-B\&B provide the same
type of additive $\varepsilon$-global guarantee on a compact scale domain, and
both have polynomial dependence on the data dimensions when $L_f$ is fixed
but exponential worst-case dependence on $L_f$. Their certificates arise
differently. NL-GRID combines a uniform Lipschitz bound with the covering
radius of a predetermined grid, whereas NL-B\&B derives adaptive upper bounds
from McCormick--ECP relaxations of entire scale boxes and does not require an
explicit global Lipschitz constant. The worst case offers no dominance:
NL-B\&B may need to refine essentially every box and its node relaxation is
larger than a fixed-scale NL-ECP, while NL-GRID always pays for every grid
point. In practice, NL-B\&B is typically preferable when bounds are strong
because incumbent solutions, bound tightening, and pruning concentrate work
near competitive scale vectors and provide a valid optimality gap throughout
the run. NL-GRID remains attractive when there are very few free scales, the
requested tolerance is coarse, or substantial parallel hardware is available,
since its fixed-scale ECPs are simple and embarrassingly parallel. Thus the
grid method is a transparent certified baseline, while the advanced B\&B
method is designed to achieve the same guarantee with substantially less
unnecessary search on difficult estimation instances.

\section{Experimental Results}\label{sec:results}
Our experiments address three questions. First, how does the cost of certified
estimation scale with the dimensions of the problem---the number of
observations, the number of utility coefficients, and the number of
nests---and how does NL-B\&B compare with the direct scale-grid method of
Section~\ref{sec:scale-grid}? We answer this on synthetic instances, where each
dimension can be varied independently while the others are held fixed. Second,
can NL-B\&B reliably solve the nonconvex nested-logit likelihood on real data,
for which local optimization offers no global guarantee? Third, does the
additional substitution structure of nested logit improve fit and out-of-sample
prediction enough to justify its computational cost? We answer these questions
in that order. The appendix gives the data construction, nest-selection
procedure, implementation settings, and additional scaling and calibration
results.

The theory and optimization formulations use average log-likelihood.
For the plots of likelihood trajectories and the likelihood-ratio statistics,
we report summed values; these equal $N$ times their average counterparts.
Thus both bounds are multiplied by $N$ together and all relative gaps and
percentage improvements are unchanged by the normalization. We write $N_{\mathrm{full}}$ for the number of observations in a cleaned
full dataset and $N$ for the training-subsample size. Thus a sampling fraction
$a\in\{0.2,0.5,0.8\}$ gives $N\approx aN_{\mathrm{full}}$ after integer
rounding. 
Algorithm~\ref{alg:global-NL} states the additive version of the certified
objective-gap stopping rule. In the experiments the same gap is compared with
a relative target, $UB-LB\leq\varepsilon_{\mathrm{rel}}|UB|$, rather than a
fixed additive target; interval widths alone never establish certification.
We report the dimensionless relative certificate
$g_{\mathrm{rel}}=(UB-LB)/|UB|$ when $UB\neq0$, and use
$\varepsilon_{\mathrm{rel}}>0$ for its target tolerance. A relative gap of
$10^{-2}$ is $1\%$; it differs from the additive tolerance $\varepsilon$ in
the theoretical statements. If $UB=0$, the absolute gap $UB-LB$ is used.
All models are estimated by maximum likelihood under a common 1-hour limit.
NL-B\&B solves its conic relaxations and fixed-dissimilarity subproblems with
\textsc{Mosek 11.2} \citep{MOSEKcookbook}; the MNL predictive benchmark is the special case in
Corollary~\ref{cor:mnl-ecp}. The bound-constrained limited-memory Broyden--Fletcher--Goldfarb--Shanno
method (L-BFGS-B),
sequential least-squares programming (SLSQP), and trust-region constrained
(trust-constr) minimizers in \textsc{SciPy 1.16.3} provide the local nested-logit
baselines, each with five random restarts. All experiments were conducted on an Apple Mac mini (M4, 10-core
central processing unit (CPU), 16 gigabytes of random-access memory (RAM)) running macOS, with all models implemented in Python 3.13. For every nested-logit instance,
all model selection---including enumeration and selection of the nest
partition---uses the training observations only. Held-out observations are
used once, after estimation and model selection, for predictive evaluation.
Full data sources, preprocessing, parameter
settings, nest selection, and additional calibration and scaling results appear
in the appendix under ``Experimental Details and Additional Results.''

\subsection{Synthetic Instances: Problem Dimensions and Method
Comparison}\label{sec:synthetic}
Synthetic instances let us move one problem dimension at a time. An instance is
indexed by the previously defined counts $(N,J,p,L)$: observations,
alternatives, utility coefficients, and nests. Attributes $\ba_{nj}\in\mathbb{R}^{p}$ are drawn i.i.d.\ uniform on
$[-3,3]^{p}$. Each individual receives a choice set by drawing a size uniformly
from $\{2,\dots,J\}$ and then a uniformly random subset of that size, so
choice-set width is heterogeneous across the sample; the minimum of two
alternatives ensures every individual is informative, since a singleton choice
set is chosen with probability one. The observed choice is drawn uniformly from
the choice set. The nest structure is built by drawing one random ordering of
the alternatives and partitioning it into $L$ contiguous blocks whose sizes
differ by at most one, so that every alternative belongs to exactly one nest
(for $J=10$: $\{5,5\}$, $\{4,3,3\}$ and $\{3,3,2,2\}$).

Crucially, the attributes, choice sets and observed choices are generated once
per instance and reused for every nest count, as is the alternative ordering.
Instances sharing an index therefore differ \emph{only} in where the nest
boundaries fall, and the comparison across $L$ is paired: any difference in
solver behaviour is attributable to the nest structure alone rather than to
sampling variation. We use $J=10$ and the full factorial
$N\in\{100,200,500\}$, $p\in\{2,5,8\}$, $L\in\{2,3,4\}$ with ten instances per
cell, giving 270 instances, solved under the same $\varepsilon_{\mathrm{rel}}=10^{-2}$
tolerance and one-hour limit as the rest of the study.

\begin{figure}[t]
    \centering
    \includegraphics[width=\linewidth]{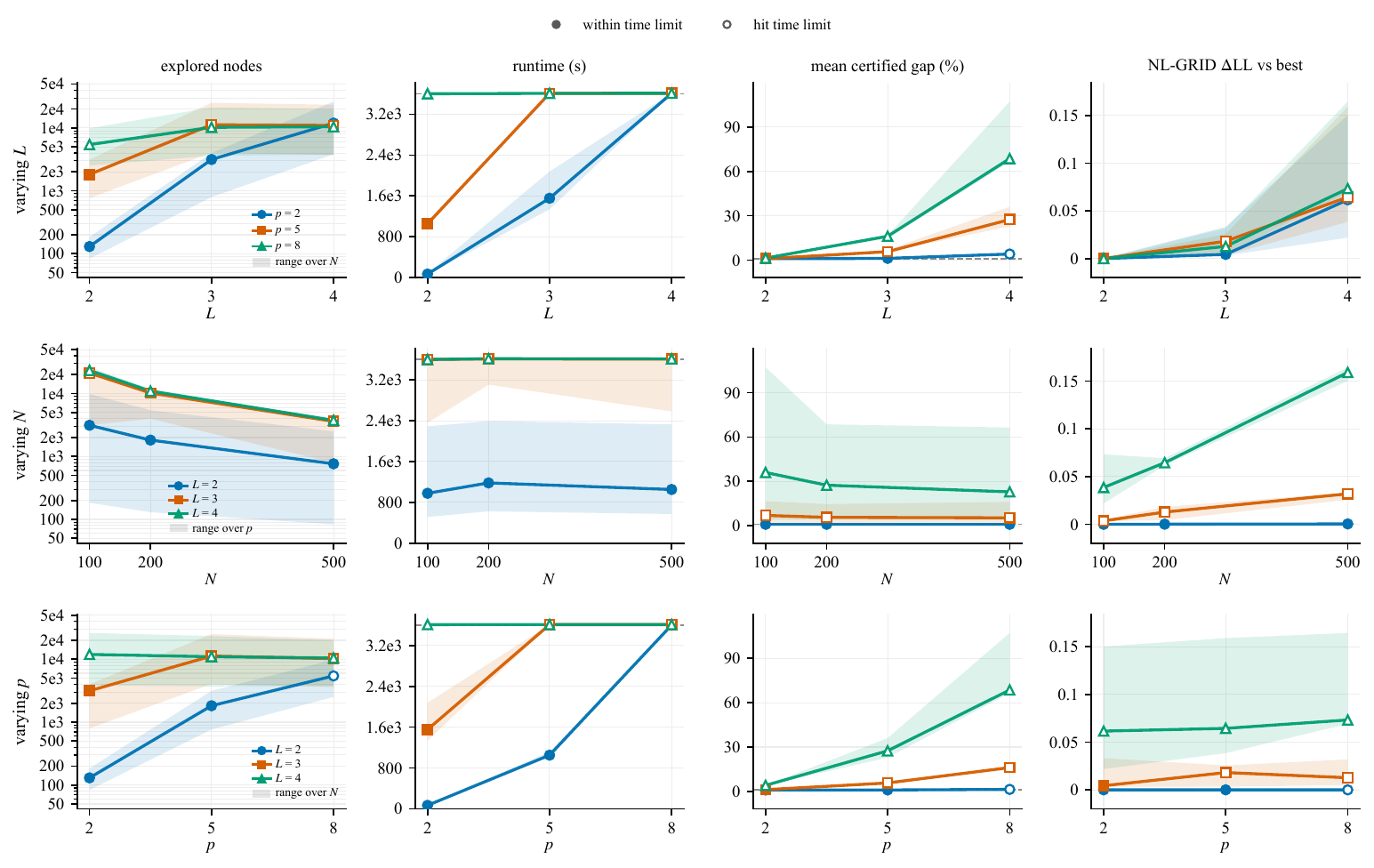}
    \caption{Scaling of certified estimation across 270 synthetic instances:
rows vary the nest count $L$, sample size $N$, and utility dimension $p$;
columns report search effort, runtime, relative gap, and NL-GRID likelihood
shortfall. The relative tolerance is $1\%$ and the time limit is one hour.}
    \label{fig:ablation-summary}
\end{figure}

Figure~\ref{fig:ablation-summary} summarizes the factorial study. Within each
row, colours identify the dimension named in the legend, lines show medians
across the remaining dimension, and shaded bands show its range. In the runtime
column, the band instead extends halfway from the median to each extreme to
keep the variation legible. Filled markers identify cells certified within the
limit; hollow markers identify cells stopped by the clock. For the latter,
reported runtimes and explored-node counts are lower bounds on the effort
needed for certification. In the plot labels, LL means log-likelihood and
$\Delta$LL is the deficit relative to the better method on the same instance.

Figure~\ref{fig:ablation-summary} separates the two factors of the search cost.
Total effort is the product of the number of nodes explored and the cost of each
node, and the two respond to different dimensions. Along the $N$ row, the median
per-node cost rises from $0.18$ to $1.03$ seconds---close to linear in the
sample size, $\approx N^{1.07}$, as the conic subproblem grows---while the
certified fraction is flat ($35.6\%$, $35.6\%$, $31.1\%$) and the mean gap does
not deteriorate ($19.5\%\to13.3\%$). More data makes each node dearer without
making the tree harder, exactly as Proposition~\ref{prop:nl-bnb-complexity-main}
predicts.

The $L$ and $p$ rows behave differently. Median node counts rise from $1{,}678$
to $10{,}871$ in $L$ and from $2{,}426$ to $9{,}320$ in $p$, but the decisive
column is the reported gap: its mean grows from $1.12\%$ to $37.87\%$
as $L$ increases from $2$ to $4$, and from $2.09\%$ to $32.64\%$ as $p$
increases from $2$ to $8$, with the certified fraction collapsing to $2.2\%$ and $8.9\%$. The
mechanisms differ. The nest count is the global-search dimension and multiplies
the tree directly; because instances sharing an index differ only in nest
structure, we can measure this without confounding, and the paired median node
ratio is $\times13.2$ from $L=2$ to $L=3$ and $\times17.2$ from $L=3$ to $L=4$.
The coefficient dimension enters only through $\bbt$, which is handled by
bound tightening rather than branching, so it leaves the per-node cost flat
($0.52$, $0.34$, $0.37$ seconds) and instead degrades the bounds. The
exponential dependence on the number of nests predicted in
Section~\ref{sec:scale-grid} is thus visible directly, and it is the only
dimension that exhibits it.

One caveat governs how these curves should be read. Instances that do not
certify within the hour are stopped by the clock, so their runtime reports the
budget and their node count only what fitted inside it. Hollow markers identify
these cells. In particular the downward slope of the node counts along the $N$
row is partly survivorship: at $N=500$ fewer instances certify, so the median is
taken over easier survivors. The per-node cost and the certified fraction are
unaffected by this and carry the scaling claim.

\begin{figure}[t]
    \centering
    \includegraphics[width=\linewidth]{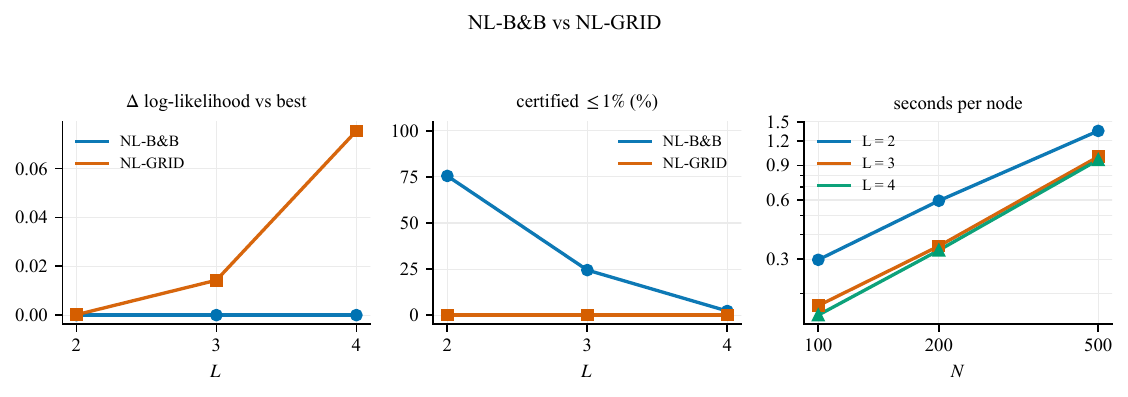}
    \caption{NL-B\&B and NL-GRID on 270 synthetic instances: likelihood deficit, certification rate, and per-node cost.}
    \label{fig:ablation-methods}
\end{figure}

Figure~\ref{fig:ablation-methods} compares three aspects of the methods.
The left panel reports the median log-likelihood deficit relative to the better
method, grouped by nest count. The centre panel reports the fraction certified
within $\varepsilon_{\mathrm{rel}}=10^{-2}$ under the one-hour limit.
The right panel shows median seconds per B\&B node against sample size on
log--log axes, with one line per nest count.

Figure~\ref{fig:ablation-methods} compares the two certified methods directly.
NL-B\&B attains the better log-likelihood on 258 of the 270 instances, and
NL-GRID's deficit widens with the nest count
($2\times10^{-4}$, $1.4\times10^{-2}$, $7.5\times10^{-2}$), since a mesh of
fixed width resolves an $L$-dimensional box ever more coarsely as $L$ grows.
The certification panel shows a larger difference: NL-B\&B certifies
$34\%$ of the instances within the one-hour limit, whereas NL-GRID certifies
none. The grid certificate includes the additive discretization allowance
$G_\lambda h/2$, where $G_\lambda$ bounds the sensitivity of the average
log-likelihood, as in Proposition~\ref{prop:grid-error}. Together with the
conic error allowance $\eta$, it gives
$UB_h-LB_h\leq G_\lambda h/2+\eta$. The experimental target is relative:
certification requires $(UB_h-LB_h)/|UB_h|\leq10^{-2}$, so the mesh needed
to attain it depends on the likelihood magnitude as well as the sensitivity
bound. A sensitivity bound computed for the summed log-likelihood must be
divided by $N$ before use with this average-objective certificate.
Although the grid-point subproblems are solved accurately, the resulting
global bounds do not meet the relative stopping criterion within the budget.
Thus NL-GRID produces competitive feasible estimates, with incumbents tracking
NL-B\&B's, but does not attain the requested certification accuracy in these
experiments. The rightmost
panel makes the complementary point about NL-B\&B: seconds per node grows as
$N^{1.07}$ and the three nest-count curves lie almost on top of one another, so
the sample size sets the cost of a node while the nest count sets how many nodes
are required. This is precisely the separation that motivates bounding whole
boxes of scale vectors rather than enumerating them.

\subsection{Real Discrete-Choice Datasets}\label{sec:certified}
We evaluate the estimators on eight established discrete-choice datasets:
Optima, SwissMetro, Cracker, ModeCanada, Heating, Fishing, Yogurt, and Mode.
They cover transportation, product choice, and household decisions, with
heterogeneous menu sizes and attribute dimensions. Each cleaned dataset is represented
in a common conditional-choice format. We generate ten random subsamples at
each of $20\%$, $50\%$, and $80\%$ of the observations, giving 30 estimation
instances per dataset; the unused observations provide held-out choices for
predictive evaluation.

Table~\ref{tab:bnb-certification} aggregates ten instances at each of the
$0.2N_{\mathrm{full}}$, $0.5N_{\mathrm{full}}$, and $0.8N_{\mathrm{full}}$
sampling levels. The improvement over a local solver is the percentage
shortfall $100(\mathcal L_{\mathrm{ref}}-\mathcal L_{\mathrm{loc}})/|\mathcal L_{\mathrm{ref}}|$,
where $\mathcal L_{\mathrm{ref}}$ is the certified incumbent's training
log-likelihood and $\mathcal L_{\mathrm{loc}}$ is a local restart's value.
We take the median of each solver's five restarts per instance, then average
over the 30 instances. A dash indicates an improvement below $0.01\%$.
For runtime presentation, each training-fraction mean is capped at
$3{,}600$ seconds when system delays cause an overrun. The dataset-level time
shown here is the arithmetic mean of the three displayed fraction means in
Table~\ref{tab:nl-cert}, rounded to seconds.

\begin{table}[t]
\centering
\footnotesize
\setlength{\tabcolsep}{3pt}
\begin{tabular}{@{}lrrrrr@{}}
\toprule
& \multicolumn{2}{c}{Search cost} & \multicolumn{3}{c}{Improvement over (\%)} \\
\cmidrule(lr){2-3}\cmidrule(l){4-6}
Dataset & Nodes & Time (s) & L-BFGS-B & SLSQP & trust-constr \\
\midrule
Mode & 198 & 130 & 0.07 & -- & 0.07 \\
Heating & 277 & 232 & 1.59 & 0.45 & 0.39 \\
Fishing & 38 & 70 & -- & -- & 0.15 \\
Optima & 556 & 1{,}120 & 0.20 & -- & -- \\
Yogurt & 26 & 94 & 0.70 & 0.01 & -- \\
ModeCanada & 146 & 928 & 0.07 & 40.98 & 1.41 \\
Cracker & 59 & 430 & 0.14 & -- & -- \\
SwissMetro & 217 & 2{,}213 & -- & -- & 0.05 \\
\bottomrule
\end{tabular}
\caption{Search cost of NL-B\&B and likelihood improvement over local solvers across 30 subsamples per dataset ($1\%$ relative tolerance).}
\label{tab:bnb-certification}
\end{table}

NL-B\&B certifies all 240 instances to the prescribed relative tolerance,
with a mean of approximately $652$ seconds across the displayed
dataset-level runtime summaries. The maximum displayed subgroup runtime is
one hour. Runtime presentation follows the same convention in both tables:
each training-fraction mean is capped at $3{,}600$ seconds when system delays
cause an overrun, and the main table averages the three displayed fraction
means. These presentation summaries do not alter objective values or gaps. The search itself stays small---between 24 and 629 nodes on average per sub-sample---and, decisively, does not grow with the sample: node counts are essentially flat across $0.2N_{\mathrm{full}}$, $0.5N_{\mathrm{full}}$ and $0.8N_{\mathrm{full}}$, so sample size enters only through the cost of each conic relaxation. What that certificate buys is a strictly better estimate than local solvers return. Averaged over the 240 instances, the certified reference exceeds a typical restart of L-BFGS-B by $0.35\%$, of trust-constr by $0.26\%$ and of SLSQP by $5.18\%$ in log-likelihood, and the improvement exceeds $0.1\%$ on 96, 27 and 25 instances respectively, reaching $95.2\%$ on ModeCanada---where SLSQP's median restart is $23$--$54\%$ below the certified reference, so it is the majority of its runs that fail, not an unlucky one. Restarting narrows the loss but does not close it: even the best of five restarts still falls more than $0.1\%$ short on 42, 6 and 2 of the 240 instances---up to $4.31\%$, on a Heating instance whose fitted reference lies at the boundary of the parameter space---and on two instances all three solvers fall short at once. None of it is flagged by the solvers' own termination status. NL-B\&B therefore certifies an estimate to the requested tolerance or returns a bound on its remaining objective error.

\begin{figure}[t]
    \centering
    \includegraphics[width=0.7\linewidth]{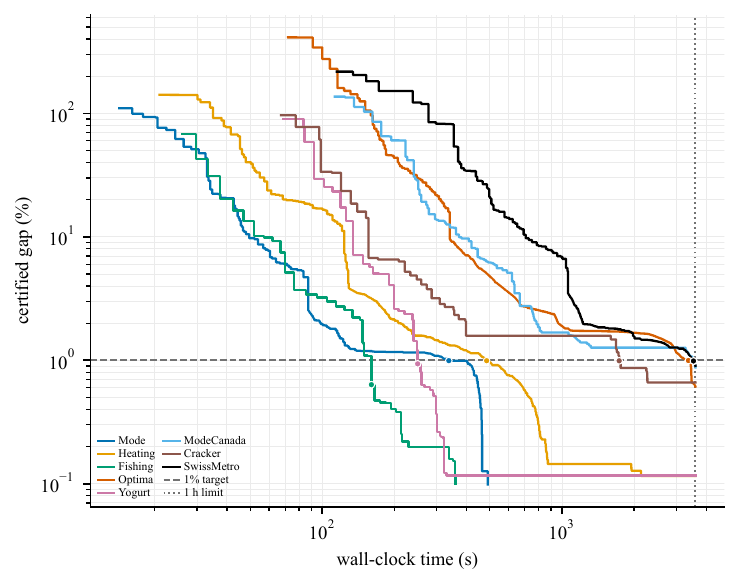}
    \caption{Certified relative optimality gaps over time on the eight full datasets (log--log axes).}
    \label{fig:convergence}
\end{figure}

The convergence experiment in Figure~\ref{fig:convergence} uses one run on
the full sample of each dataset at the tighter tolerance
$\varepsilon_{\mathrm{rel}}=10^{-3}$. It plots the certified relative gap
$(UB-LB)/|UB|$ against wall-clock time. These full-sample runs are separate
from the 240 subsample instances in Table~\ref{tab:bnb-certification}.

Figure~\ref{fig:convergence} shows how the certificate tightens over time. The bound starts loose---the reported relative gap is $68$--$413\%$ at the root---and collapses by two orders of magnitude within the first few hundred seconds, after which 6 of 8 curves flatten against the $1\%$ line. What the flattening costs differs sharply by dataset: Yogurt and Heating reach $1\%$ after only $7\%$ and $13\%$ of their eventual runtime, whereas SwissMetro and ModeCanada need $98\%$ of it, and Optima $94\%$. Certification is therefore front-loaded in accuracy but back-loaded in effort---the first $1\%$ of the gap is cheap, the last fraction of a percent is not---the incumbent often stabilizes before its global accuracy has been certified.

\subsection{Predictive Value of Nested Logit}\label{sec:predictive}
\begin{figure}[t]
    \centering
    \includegraphics[width=0.7\linewidth]{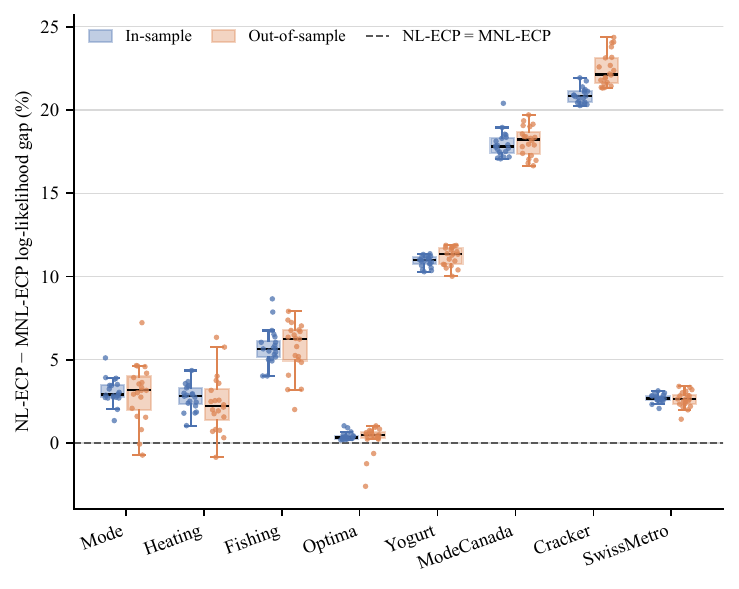}
    \caption{Relative log-likelihood improvement of nested logit over MNL on training and held-out samples.}
    \label{fig:NLvsMNL}
\end{figure}

In Figure~\ref{fig:NLvsMNL}, boxes summarize the resamples and points show
individual instances. NL and MNL denote the jointly fitted nested-logit and
multinomial-logit models; the zero line indicates equal evaluation
log-likelihoods.

Figure~\ref{fig:NLvsMNL} compares nested logit with MNL using
$(\mathcal L_{\mathrm{NL}}^{\mathrm{eval}}-\mathcal L_{\mathrm{MNL}}^{\mathrm{eval}})/|\mathcal L_{\mathrm{MNL}}^{\mathrm{eval}}|$ over ten
$0.5N_{\mathrm{full}}$ and $0.8N_{\mathrm{full}}$ resamples per dataset. Here $\mathcal L_{\mathrm{NL}}^{\mathrm{eval}}$ and
$\mathcal L_{\mathrm{MNL}}^{\mathrm{eval}}$ are the summed log-likelihoods of
the fitted models on the same evaluation sample (training or held-out, as
indicated). These evaluation sums are distinct from the single-observation
log-likelihood $\ell(\theta;X)$ in the statistical analysis. Nested logit improves in-sample fit on
all 160 instances, with a mean gain of $8.1\%$ (from $0.4\%$ on Optima to more
than $21\%$ on Cracker). More importantly, the improvement transfers to unseen
choices: nested logit achieves the higher held-out likelihood on 154 of 160
instances ($96\%$), and its mean out-of-sample gain of $8.2\%$ closely matches
the in-sample gain. The few reversals occur only where the nesting signal is
weak and the difference is near zero. Hence the extra scale parameters capture
stable substitution patterns rather than merely increasing in-sample
flexibility. Together with Table~\ref{tab:bnb-certification}, this result shows
the value of our approach: NL-B\&B makes the statistically superior nested-logit
model usable without sacrificing a verifiable global-accuracy guarantee.

\section{Conclusion}
\label{sec:conclusion}
\label{sec:concl}
We provide the first globally certified framework for joint maximum-likelihood
estimation of nested logit. The formulation exploits a sharp separation in the
problem structure: utility estimation is an exponential-cone program once nest
dissimilarities are fixed, while all nonconvex search is confined to those
dissimilarities. This yields valid anytime bounds, finite $\varepsilon$-global
convergence on a compact domain, and an additive PTAS when the number of free
nests is fixed. The guarantee complements the OR literature on assortment,
pricing, revenue management, and facility location, which typically optimizes
decisions conditional on estimated choice parameters. The computational study shows why the estimation certificate matters. Across
240 instances, NL-B\&B certifies every reported estimate, whereas multistart
local solvers sometimes terminate normally at inferior likelihoods. The
certificate quantifies the remaining estimation error in likelihood value,
providing a verifiable measure of numerical solution quality for the demand
parameters supplied to downstream optimization models.

The method assumes a prespecified nesting structure and a compact parameter
domain, and its worst-case dependence on the number of free nests is
exponential. Difficult instances stopped early retain a nonzero but valid gap.
Promising extensions include joint nest-structure selection, cross-nested and
mixed-logit estimation, parallel node processing, and decision-aware guarantees
that propagate estimation and sampling uncertainty into assortment, pricing,
capacity, and location decisions.

%\section*{Acknowledgments}

\bibliographystyle{informs2014}
\bibliography{references}

\clearpage
\appendix
\begin{center}
{\Large\bfseries Appendix: Proofs and Supplementary Material}
\end{center}
\section*{Appendix Outline}
The appendix is organized as follows. First, ``Proofs and Extensions'' proves
the fixed-dissimilarity nested-logit ECP, its MNL corollary, and the error and
complexity guarantees for direct scale-grid estimation. Second, ``Full B\&B Details''
presents the complete lifted relaxation, bound-tightening and branching rules,
convergence and complexity arguments, and a visual workflow. The intervening
section ``Sample Efficiency and Computational Complexity'' connects statistical
consistency to the polynomial dependence of estimation cost on sample size. Finally,
``Experimental Details and Additional Results'' documents the data, model
selection, solver settings, and supplementary performance results.

\section*{Appendix: Proofs and Extensions}

\paragraph{Proof of Proposition~\ref{prop:fixed-nl-ecp}.}
\label{ec:proof:fixed-nl-ecp}
For fixed $\blambda$, introduce $z_{nl}=\log W_{nl}$ for every active pair
$l\in\cI_n$ and
$y_n=\log\sum_{l\in\cI_n}\exp(\lambda_lz_{nl})$.  Substitution into
\eqref{prop:NL-MLE} gives the objective in \eqref{prob:NL-ECP}.  The two
log-sum-exp epigraphs are equivalent to
\[
 \sum_{j\in\cN_l\cap S_n}
 e^{\bbt^\top\ba_{nj}/\lambda_l-z_{nl}}\leq1,
 \qquad
 \sum_{l\in\cI_n}e^{\lambda_lz_{nl}-y_n}\leq1.
\]
Represent each exponential term by, respectively,
$(k_{njl},1,\bbt^\top\ba_{nj}/\lambda_l-z_{nl})\in\cK_{\exp}$ and
$(h_{nl},1,\lambda_lz_{nl}-y_n)\in\cK_{\exp}$, and bound the corresponding
sums by one \citep{Chares2009,MOSEKcookbook}. This yields exactly \eqref{prob:NL-ECP}. Conversely, any feasible
ECP point implies both epigraph inequalities by summing the lower bounds on
$k_{njl}$ and $h_{nl}$; any feasible epigraph point extends to the ECP by
setting those variables equal to their exponential expressions.

It remains to verify tightness.  The objective decreases with $y_n$, so its
epigraph is tight.  Reducing any $z_{nl}$ to its within-nest log-sum-exp cannot
decrease feasibility of the upper-level epigraph.  For the chosen nest its
objective coefficient is $\lambda_{l_n}-1\leq0$, and for every other nest it
affects the objective only through the minimized $y_n$; hence an optimum exists
with every within-nest epigraph tight.  The ECP and fixed-scale nested-logit MLE
therefore have identical optimal values and recover the same $\bbt$.

\paragraph{Proof of Corollary~\ref{cor:mnl-ecp}.}
\label{ec:proof:mnl-ecp}
When $\lambda_l=1$ for all $l$, the nested-logit probability simplifies to
\[
 P_n(j\mid S_n)=
 \frac{\exp(\bbt^\top\ba_{nj})}
 {\sum_{k\in S_n}\exp(\bbt^\top\ba_{nk})},
\]
independently of the declared partition. In \eqref{prob:NL-ECP}, the
coefficients $(\lambda_{l_n}-1)$ vanish and the two nested log-sum-exp
epigraphs compose into the single epigraph
$t_n\geq\log\sum_{j\in S_n}\exp(\bbt^\top\ba_{nj})$. Writing this epigraph as
$\sum_{j\in S_n}\exp(\bbt^\top\ba_{nj}-t_n)\leq1$ and introducing $u_{nj}$
gives \eqref{prob:MNL-ECP}.

For complexity, \eqref{prob:MNL-ECP} contains $Z=\sum_n|S_n|$
constant-dimensional exponential-cone blocks. Their self-concordant barrier
parameter is of order $Z$, so a short-step primal--dual method requires
$O(\sqrt Z\log(1/\varepsilon))$ Newton steps
\citep{NesterovNemirovskii1994,BenTalNemirovski2001}. Eliminating cone-local directions
costs $O(Zp^2)$ per step, and a conservative dense factorization of the Schur
system on the $p$ shared coefficients and $N$ observation epigraphs costs
$O((p+N)^3)$. Multiplying the iteration and per-step bounds proves the stated
complexity.

\paragraph{Proof of Proposition~\ref{prop:grid-error}.}
\label{ec:proof:grid-error}
We prove the Lipschitz, approximation, and certification claims in turn.
Because every attribute vector is finite, $\cB\times\cD$ is compact,
and every dissimilarity is bounded away from zero, the nested-logit
log-likelihood is continuously differentiable in the free dissimilarities on
an open neighborhood of $\cB\times\cD$. Consequently, its scale gradient is
bounded there, the supremum in Proposition~\ref{prop:grid-error} is finite,
and any $G_\lambda$ dominating that supremum is valid. For any
$\bbt\in\cB$ and $\blambda,\blambda'\in\cD$, the line segment joining the two
scale vectors lies in the box $\cD$. The multivariate mean-value inequality
and duality of the one- and infinity-norms give
\begin{equation}
 \left|\cL^{\NL}(\bbt,\blambda)
       -\cL^{\NL}(\bbt,\blambda')\right|
 \leq G_\lambda\|\blambda-\blambda'\|_\infty.
 \label{eq:grid-uniform-lipschitz}
\end{equation}

We next show that maximizing over $\bbt$ preserves this uniform Lipschitz
constant. Compactness and continuity ensure that the maximum defining
$V(\blambda)$ is attained. Let $\bbt_{\blambda}$ be a maximizer at
$\blambda$. Applying \eqref{eq:grid-uniform-lipschitz} at this coefficient
vector yields
\begin{align*}
 V(\blambda)-V(\blambda')
 &=\cL^{\NL}(\bbt_{\blambda},\blambda)
   -\max_{\bbt\in\cB}\cL^{\NL}(\bbt,\blambda')\\
 &\leq \cL^{\NL}(\bbt_{\blambda},\blambda)
   -\cL^{\NL}(\bbt_{\blambda},\blambda')\\
 &\leq G_\lambda\|\blambda-\blambda'\|_\infty.
\end{align*}
Interchanging $\blambda$ and $\blambda'$ proves
\begin{equation}
 |V(\blambda)-V(\blambda')|
 \leq G_\lambda\|\blambda-\blambda'\|_\infty.
 \label{eq:profile-lipschitz}
\end{equation}

For coordinate $l$, the construction in \eqref{eq:scale-grid} divides its
interval into $K_l$ subintervals of length
$\Delta_l/K_l\leq h$. Rounding any coordinate to its nearest endpoint thus
moves it by at most $h/2$. Hence every $\blambda\in\cD$ has a grid
representative $\pi_h(\blambda)\in\cG_h$, obtained by rounding each coordinate, satisfying
\begin{equation}
 \|\blambda-\pi_h(\blambda)\|_\infty\leq h/2.
 \label{eq:grid-covering-radius}
\end{equation}
Let $\blambda^\star$ maximize $V$ on $\cD$ and write
$\blambda^g=\pi_h(\blambda^\star)$. By
\eqref{eq:profile-lipschitz}--\eqref{eq:grid-covering-radius},
\begin{equation}
 V^\star\leq V(\blambda^g)+G_\lambda h/2.
 \label{eq:grid-optimum-neighbor}
\end{equation}

At each grid point, primal--dual feasibility and a duality gap of at most
$\eta$ imply
\begin{equation}
 \underline V_\eta(\blambda)\leq V(\blambda)
 \leq\overline V_\eta(\blambda),
 \qquad
 \overline V_\eta(\blambda)-\underline V_\eta(\blambda)\leq\eta.
 \label{eq:grid-conic-certificate}
\end{equation}
Using the definition of $LB_h$, followed by
\eqref{eq:grid-conic-certificate} and \eqref{eq:grid-optimum-neighbor}, gives
\[
 LB_h\geq \underline V_\eta(\blambda^g)
 \geq V(\blambda^g)-\eta
 \geq V^\star-G_\lambda h/2-\eta.
\]
Because every primal grid solution is feasible for the original bounded
NL-MLE, $LB_h\leq V^\star$, proving the stated lower-bound error.

For any $\blambda\in\cD$, apply the same covering argument and then the dual
bound in \eqref{eq:grid-conic-certificate}:
\[
 V(\blambda)
 \leq V(\pi_h(\blambda))+G_\lambda h/2
 \leq \max_{\bmu\in\cG_h}\overline V_\eta(\bmu)+G_\lambda h/2
 =UB_h.
\]
Taking the maximum over $\blambda$ proves that $UB_h$ is globally valid.
Moreover, \eqref{eq:grid-conic-certificate} implies
\[
 UB_h-LB_h\leq G_\lambda h/2+\eta,
\]
so the reported interval is itself a certificate of the claimed accuracy.
For $G_\lambda>0$, setting $h=\varepsilon/G_\lambda$ and
$\eta=\varepsilon/2$ makes this gap at most $\varepsilon$. If
$G_\lambda=0$, the profile likelihood is constant on $\cD$, so one arbitrary
scale vector suffices and only the conic tolerance remains.

\paragraph{Proof of Proposition~\ref{prop:grid-complexity}.}
\label{ec:proof:grid-complexity}
Coordinate $l$ contributes exactly $K_l+1$ values, and hence
\[
 |\cG_h|=\prod_{l\in\cF}
 \left(1+\left\lceil\frac{\Delta_l}{h}\right\rceil\right).
\]
The fixed-scale formulation contains $Z+\Lambda$ constant-dimensional
exponential-cone blocks. The standard short-step iteration bound is
$O(\sqrt{Z+\Lambda}\log(1/\eta))$; eliminating cone-local directions and
factoring a conservative dense Schur system costs
$O(Zp^2+(p+N+\Lambda)^3)$ per iteration
\citep{NesterovNemirovskii1994,BenTalNemirovski2001}. Multiplying this
per-grid-point cost by $|\cG_h|$ proves
\eqref{eq:grid-complexity-general}. Substituting
$h=\varepsilon/G_\lambda$ and $\eta=\varepsilon/2$ gives
\eqref{eq:grid-complexity}. For fixed $L_f$, the product contains a fixed
number of factors, each polynomial in the data dimensions and
$1/\varepsilon$ (under the stated compact-domain derivative bound). If $L_f$
is part of the input, the Cartesian-product factor is exponential in $L_f$.
Finally, the grid-point ECPs share no decision variables or constraints, so
they may be dispatched in parallel. Such parallelism can reduce elapsed time
but does not alter the sequential total-work bound proved above.

\section*{Appendix: Sample Efficiency and Computational Complexity}
\label{app:sample-complexity}
This section connects the statistical benefit of collecting more observations
to the cost of globally solving the resulting nested-logit empirical
likelihood. Fix a target population error $\varepsilon_{\mathrm{stat}}>0$ and failure probability $\alpha\in(0,1)$. Let $X_1,\ldots,X_N$ be i.i.d. choice observations (including the
menu, attributes, and selected alternative), let
$\theta=(\bbt,\blambda_{\cF})\in\Theta$ collect the utility coefficients and the
$L_f$ free nonsingleton nest scales, and let $\ell(\theta;X_n)$ be one
observation's log-likelihood. Define
\[
 \begin{aligned}
 \widehat Q_N(\theta)
   &=\frac1N\sum_{n=1}^N\ell(\theta;X_n),\\
 Q(\theta)
   &=\mathrm{E}[\ell(\theta;X)].
 \end{aligned}
\]
For the observed sample,
$\widehat Q_N(\theta)=\cL^{\NL}(\bbt,\blambda)$. This is the same
average objective used by NL-ECP, NL-GRID, and NL-B\&B.

\paragraph{Regularity conditions.}
We use standard conditions for maximum-likelihood extremum estimators
\citep{NeweMcFa94}. Let $\Theta$ be a compact subset of
$\mathbf R^{p+L_f}$ with Euclidean diameter $R$, and let $\theta_0$ be any maximizer
of $Q$. Let $G>0$ be a parameter-Lipschitz bound and $B>0$ an
observation-range bound. Assume that, uniformly over observations $x$ in the support of $X$,
\[
 |\ell(\theta;x)-\ell(\theta';x)|
 \leq G\|\theta-\theta'\|_2
 \quad\text{for all }\theta,\theta'\in\Theta,
\]
and, uniformly over $\theta\in\Theta$, that
\[
 \operatorname*{ess\,sup}_{x}\ell(\theta;x)
 -\operatorname*{ess\,inf}_{x}\ell(\theta;x)\leq B.
\]
Here $\operatorname*{ess\,sup}$ and $\operatorname*{ess\,inf}$ are the
supremum and infimum after ignoring observation events of probability zero.
For nested logit, these conditions follow, for example, from bounded
attributes, bounded coefficient and scale domains, bounded menus, and choice
probabilities uniformly bounded away from zero.

For the computational corollary we use the same uniform relaxation-error
condition as Algorithm~\ref{alg:global-NL}: widths at most $\delta$ and
endpoint ratios at most $1+\delta$, with combined conic allowance $\eta$,
give certified average-objective error at most $H\delta+\eta$. The constant
$H$ is independent of $N$ on the declared bounded domains. This hypothesis
is needed for the stated B\&B approximation and complexity guarantees, but
not for the statistical statement in Proposition~\ref{prop:sample-computation}.
An error of $\varepsilon$ in this average objective equals $N\varepsilon$
in the corresponding summed objective.

\paragraph{From population accuracy to sample size.}
Let $N^\star(\varepsilon_{\mathrm{stat}},\alpha)$ denote a sufficient sample
size for the prescribed population accuracy and failure probability. A
covering-number argument and Hoeffding's inequality imply that it suffices to use
\[
 N^\star(\varepsilon_{\mathrm{stat}},\alpha)
 =O\!\left(\frac{B^2}{\varepsilon_{\mathrm{stat}}^2}
 \left[(p+L_f)\log\!\left(1+\frac{GR}{\varepsilon_{\mathrm{stat}}}\right)
       +\log\frac1\alpha\right]\right)
\]
observations to make the uniform statistical error at most $\varepsilon_{\mathrm{stat}}/3$
with probability at least $1-\alpha$ (constant factors are absorbed by $O(\cdot)$).
If an algorithm then returns $\widehat\theta_N$ whose average empirical
log-likelihood is within $\varepsilon_{\mathrm{stat}}/3$ of the empirical maximum, the usual
three-term decomposition gives
\[
 \Pr\!\left\{Q(\theta_0)-
 Q(\widehat\theta_N)\leq\varepsilon_{\mathrm{stat}}\right\}\geq1-\alpha.
\]
\paragraph{Proof of Proposition~\ref{prop:sample-computation}.}
\label{ec:proof:sample-computation}
Cover $\Theta$ by a Euclidean net of radius proportional to
$\varepsilon_{\mathrm{stat}}/G$. Its cardinality is at most
$(1+c_0GR/\varepsilon_{\mathrm{stat}})^{p+L_f}$ for a universal constant $c_0$.
Hoeffding's inequality and a union bound over this net, followed by the
$G$-Lipschitz extension from the net to all of $\Theta$, show that the stated
sample size ensures, with probability at least $1-\alpha$,
\[
 \Delta_N:=\sup_{\theta\in\Theta}
 |\widehat Q_N(\theta)-Q(\theta)|\leq\varepsilon_{\mathrm{stat}}/3.
\]
On this event, empirical $\varepsilon_{\mathrm{stat}}/3$-optimality gives
\[
 \Gamma_N:=\sup_{\theta\in\Theta}\widehat Q_N(\theta)
 -\widehat Q_N(\widehat\theta_N)\leq\varepsilon_{\mathrm{stat}}/3.
\]
Therefore,
\[
 Q(\theta_0)-Q(\widehat\theta_N)
 \leq2\Delta_N+\Gamma_N\leq\varepsilon_{\mathrm{stat}}.
\]
This combines statistical and optimization error in excess expected
log-likelihood. It also shows that there is little benefit in solving the
empirical problem much more accurately than its $O(N^{-1/2})$ statistical
error. A parameter-distance or choice-probability guarantee would require
additional identification or curvature assumptions.

\paragraph{Proof of Corollary~\ref{cor:sample-computation}.}
\label{ec:proof:sample-computation-cor}
The arithmetic bounds follow by substituting the required sample size and an
$O(\varepsilon_{\mathrm{stat}})$ average-objective tolerance into
Proposition~\ref{prop:nl-bnb-complexity-main} and using the quantitative
relaxation condition above.

\paragraph{End-to-end complexity for nested logit.}
For nested logit, bounded menus give $Z=O(N)$ and $\Lambda=O(NL)$. Let
$N=N^\star(\varepsilon_{\mathrm{stat}},\alpha)$. Under the quantitative relaxation
condition above, choose $\delta=\varepsilon_{\mathrm{stat}}/(6H)$ and
$\eta=\varepsilon_{\mathrm{stat}}/6$. The resulting certified optimization error in average
log-likelihood is at most $\varepsilon_{\mathrm{stat}}/3$. Use $D=O(p+Lp+L+N^\star(1+L))$ for the global linear-system dimension.
Let $\mathcal C_{\mathrm{NL}}$ denote the resulting arithmetic-operation count.
Proposition~\ref{prop:nl-bnb-complexity-main} then yields
\[
 \begin{aligned}
 \mathcal C_{\mathrm{NL}}(\varepsilon_{\mathrm{stat}},\alpha)
 =O\!\Bigg(&\left(1+\frac{\kappa H}{\varepsilon_{\mathrm{stat}}}\right)^{L_f}
 \sqrt{N^\star(1+L)}\log\frac1\varepsilon_{\mathrm{stat}}\\[-1mm]
 &\cdot\Big[N^\star p^2+D^3\Big]\Bigg).
 \end{aligned}
\]
Thus, with probability at least $1-\alpha$,
$Q(\theta_0)-Q(\widehat\theta_N)\leq\varepsilon_{\mathrm{stat}}$.
The required computation is polynomial in $N^\star$ and
$1/\varepsilon_{\mathrm{stat}}$ for every fixed $L_f$. The
important distinction is that increasing $N$ enlarges only the convex ECPs at
each node; it does not create a new branching coordinate. The exponential
factor is controlled by the number of behavioral groups (nests), not by the
number of observations. Consequently, when the nest structure is small, the
estimator can exploit increasing samples for consistent population prediction
without incurring exponential growth in sample size.

\section*{Appendix: Full B\&B Details}
\label{app:bnb-details}
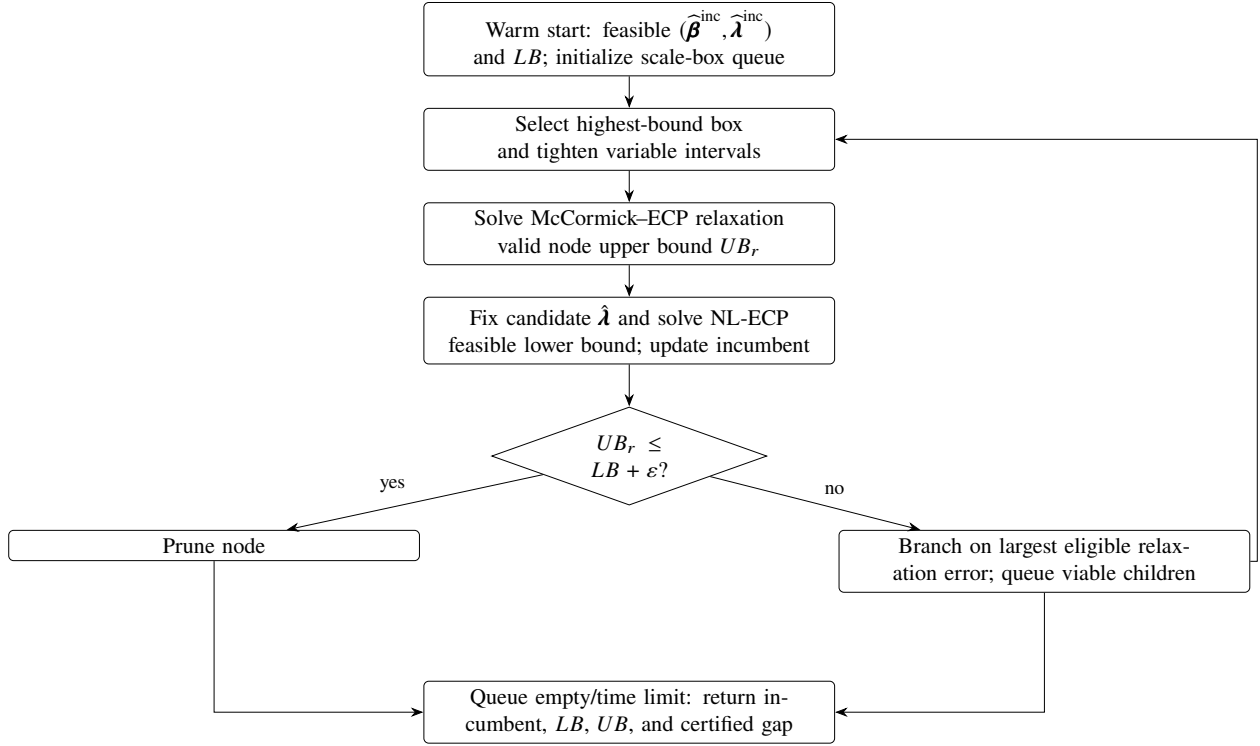
\begin{figure*}[t]
\centering
\resizebox{\textwidth}{!}{%
\begin{tikzpicture}[
  node distance=4.5mm,
  >={Stealth[length=1.8mm]},
  box/.style={draw,rounded corners=2pt,align=center,font=\footnotesize,
              text width=5.6cm,inner ysep=3pt},
  test/.style={draw,diamond,aspect=2.8,align=center,font=\footnotesize,
               text width=1.8cm,inner ysep=1pt}]
\node[box] (init) {Warm start: feasible $(\widehat\bbt^{\mathrm{inc}},\widehat\blambda^{\mathrm{inc}})$ and $LB$; initialize scale-box queue};
\node[box,below=of init] (select) {Select highest-bound box and tighten variable intervals};
\node[box,below=of select] (relax) {Solve McCormick--ECP relaxation\\valid node upper bound $UB_r$};
\node[box,below=of relax] (inc) {Fix candidate $\hat\blambda$ and solve NL-ECP\\feasible lower bound; update incumbent};
\node[test,below=6mm of inc] (prune) {$UB_r\leq LB+\varepsilon$?};
\node[box,below left=7mm and 20mm of prune] (discard) {Prune node};
\node[box,below right=7mm and 20mm of prune] (branch) {Branch on largest eligible relaxation error; queue viable children};
\node[box,below=25mm of prune] (done) {Queue empty/time limit: return incumbent, $LB$, $UB$, and certified gap};
\draw[->] (init)--(select);
\draw[->] (select)--(relax);
\draw[->] (relax)--(inc);
\draw[->] (inc)--(prune);
\draw[->] (prune)--node[above left,font=\scriptsize]{yes}(discard);
\draw[->] (prune)--node[above right,font=\scriptsize]{no}(branch);
\draw[->] (discard.south) |- (done.west);
\draw[->] (branch.south) |- (done.east);
\draw[->] (branch.east) -- ++(1mm,0) |- (select.east);
\end{tikzpicture}%
}
\caption{Workflow of globally certified nested-logit estimation by NL-B\&B.}
\label{fig:bnb-flow}
\end{figure*}

Figure~\ref{fig:bnb-flow} summarizes the interaction between the feasible
incumbent, convex node bounds, pruning, and branching.

We now consider the joint estimation of the utility parameters \(\bbt\) and the
nested-logit scale parameters \(\blambda=(\lambda_1,\ldots,\lambda_L)\).
Recall that, for fixed \(\blambda\), the nested-logit likelihood admits an
exponential-cone reformulation and can be solved to global optimality. However,
when \(\blambda\) is also optimized, the problem becomes nonconvex because the
likelihood and the exponential-cone constraints contain products and fractional
terms involving \(\lambda_l\). In particular, the formulation contains terms of
the form
\[
    \lambda_l z_{nl},
    \qquad
    \frac{\bbt^\top \ba_{nj}}{\lambda_l},
\]
which destroy convexity when \(\lambda_l\) is treated as a decision variable.

{\color{black}\paragraph{Selecting and validating parameter bounds.}
Finite bounds are part of the stated estimation problem, rather than a hidden
numerical device. Attributes are first standardized so coefficient magnitudes
are comparable. An initial coefficient box can be obtained from substantive
marginal-utility limits or from a preliminary multi-start fit, enlarged by a
conservative margin based on coefficient magnitude and standard error. Scale
bounds use the random-utility restriction $\lambda_l\leq1$ and a small positive
lower limit chosen before examining held-out outcomes. A likelihood-level-set
bound may also be constructed from a feasible incumbent: interval propagation
and optimization-based bound tightening retain every parameter value whose
likelihood can exceed the incumbent while removing only regions certified not
to do so.

We validate the resulting domain by a systematic expansion test. After solving
the bounded problem, double every active coefficient radius and decrease each
active scale lower bound, then solve and certify again. Expansion continues
until (i) the incumbent is interior, (ii) its parameter vector and likelihood
are unchanged to the reporting tolerance, and (iii) the certified upper bound
does not reveal an improving region in the added domain. These checks show
stability with respect to the selected box, but the formal certificate always
applies to the final reported domain. If expansion causes coefficients to
diverge, the unconstrained MLE is not finite (for example under separation),
and the bounded or regularized estimand must be reported explicitly.}

To obtain a global solution method, we propose a spatial branch-and-bound
algorithm \citep{HorstTuy1996,TawarmalaniSahinidis2002} that branches on the dissimilarity parameters \(\lambda_l\). At each
node of the branch-and-bound tree, the variables \(\lambda_l\) are restricted to
a box
\[
    \lambda_l \in [\underline{\lambda}_l^r,\overline{\lambda}_l^r],
    \qquad l\in[L],
\]
where \(r\) denotes the current node. Over this box, the nonconvex terms are
replaced by convex relaxations based on McCormick envelopes
\citep{mccormick1976computability}. The resulting node
relaxation is an exponential-cone program, which provides a valid upper bound on
the optimal log-likelihood over that node. Feasible solutions, and hence lower
bounds, are obtained by fixing \(\blambda\) and solving the corresponding
convex NL-ECP in \(\bbt\).

\paragraph{Reformulation of the nonconvex terms.}

The main source of nonconvexity in the NL-ECP formulation is the term
\(\lambda_l z_{nl}\) appearing in the upper-level exponential-cone constraint
\[
    (h_{nl},1,\lambda_l z_{nl}-y_n)\in\cK_{\exp}.
\]
To linearize this product, we introduce an auxiliary variable
\[
    q_{nl} = \lambda_l z_{nl},
    \qquad n\in[N],\ l\in\cI_n.
\]
At a branch-and-bound node \(r\), let
$\underline z_{nl}^r,\overline z_{nl}^r$ be valid finite bounds for the
inclusive-value epigraph variable, for each active pair $l\in\cI_n$.
Suppose that
\[
    \lambda_l \in [\underline{\lambda}_l^r,\overline{\lambda}_l^r],
    \qquad
    z_{nl} \in [\underline{z}_{nl}^r,\overline{z}_{nl}^r].
\]
Then the equality \(q_{nl}=\lambda_l z_{nl}\) can be relaxed by the McCormick
envelope \citep{mccormick1976computability}:
\begin{align}
    q_{nl}
    &\ge
    \underline{\lambda}_l^r z_{nl}
    + \underline{z}_{nl}^r \lambda_l
    - \underline{\lambda}_l^r \underline{z}_{nl}^r,
    \label{eq:mccormick-1}
    \\
    q_{nl}
    &\ge
    \overline{\lambda}_l^r z_{nl}
    + \overline{z}_{nl}^r \lambda_l
    - \overline{\lambda}_l^r \overline{z}_{nl}^r,
    \label{eq:mccormick-2}
    \\
    q_{nl}
    &\le
    \underline{\lambda}_l^r z_{nl}
    + \overline{z}_{nl}^r \lambda_l
    - \underline{\lambda}_l^r \overline{z}_{nl}^r,
    \label{eq:mccormick-3}
    \\
    q_{nl}
    &\le
    \overline{\lambda}_l^r z_{nl}
    + \underline{z}_{nl}^r \lambda_l
    - \overline{\lambda}_l^r \underline{z}_{nl}^r.
    \label{eq:mccormick-4}
\end{align}
The exponential-cone constraint is then relaxed as
\[
    (h_{nl},1,q_{nl}-y_n)\in\cK_{\exp},
    \qquad n\in[N],\ l\in\cI_n.
\]

The second nonconvex term is
\[
    \frac{\bbt^\top \ba_{nj}}{\lambda_l}.
\]
To handle this term, we introduce the reciprocal variable
\[
    \rho_l = \frac{1}{\lambda_l},
    \qquad l\in[L],
\]
and auxiliary variables
\[
    \bgamma_l = \rho_l \bbt,
    \qquad l\in[L].
\]
Then
\[
    \frac{\bbt^\top \ba_{nj}}{\lambda_l}
    =
    \bgamma_l^\top \ba_{nj}.
\]
The equality \(\bgamma_l=\rho_l\bbt\) is bilinear. Hence, for each coefficient
\(d\in[p]\), we introduce
\[
    \gamma_{ld} = \rho_l \beta_d,
    \qquad l\in[L],\ d\in[p],
\]
and relax it by McCormick inequalities over the bounds
\[
    \rho_l\in[\underline{\rho}_l^r,\overline{\rho}_l^r],
    \qquad
    \beta_d\in[\underline{\beta}_d,\overline{\beta}_d],
\]
where $\underline\beta_d=\min_{\bbt\in\cB}\beta_d$ and
$\overline\beta_d=\max_{\bbt\in\cB}\beta_d$ are the declared coordinate
bounds, and the reciprocal endpoints are
\[
    \underline{\rho}_l^r = \frac{1}{\overline{\lambda}_l^r},
    \qquad
    \overline{\rho}_l^r = \frac{1}{\underline{\lambda}_l^r}.
\]
The McCormick envelope for
\(\gamma_{ld}=\rho_l\beta_d\) is
\begin{align}
    \gamma_{ld}
    &\ge
    \underline{\rho}_l^r \beta_d
    + \underline{\beta}_d \rho_l
    - \underline{\rho}_l^r\underline{\beta}_d,
    \\
    \gamma_{ld}
    &\ge
    \overline{\rho}_l^r \beta_d
    + \overline{\beta}_d \rho_l
    - \overline{\rho}_l^r\overline{\beta}_d,
    \\
    \gamma_{ld}
    &\le
    \underline{\rho}_l^r \beta_d
    + \overline{\beta}_d \rho_l
    - \underline{\rho}_l^r\overline{\beta}_d,
    \\
    \gamma_{ld}
    &\le
    \overline{\rho}_l^r \beta_d
    + \underline{\beta}_d \rho_l
    - \overline{\rho}_l^r\underline{\beta}_d.
\end{align}
Using \(\bgamma_l\), the within-nest exponential-cone constraint becomes
\[
    (k_{njl},1,\bgamma_l^\top\ba_{nj}-z_{nl})\in\cK_{\exp}.
\]

\paragraph{Node relaxation.}

For scalar factors $x,y$ and a product variable $w$, let
$w\in\operatorname{MC}_r(x,y)$ mean that $(x,y,w)$ satisfies the
four McCormick inequalities for $w=xy$ on the specified node-$r$ bounds,
as in \eqref{eq:mccormick-1}--\eqref{eq:mccormick-4}. The envelope includes
these bounds on $x$ and $y$. Let $\brho=(\rho_l)_{l\in[L]}$,
$\bgamma=(\gamma_{ld})_{l\in[L],d\in[p]}$, and
$\bq=(q_{nl})_{n\in[N],l\in\cI_n}$ collect the lifted variables.
We call the node relaxation NL-MC-ECP, where MC denotes McCormick.
At node $r$, it is

\begin{align}
   \max_{\substack{\bbt,\blambda,\brho,\bgamma,\\ \bz,\by,\bk,\bh,\bq}}
   \quad
   & \frac1N\sum_{n\in[N]}
   \big[
        q_{nl_n} - z_{n l_n}
   \label{prob:NL-MC-ECP}
   \tag{\sf NL-MC-ECP}
   \\
   & \qquad
        + \bgamma_{l_n}^{\top}\ba_{n j_n}
        - y_n
   \big]
   \nonumber
   \\
   \text{s.t.}\quad
   &
   \sum_{j\in\cN_l\cap S_n} k_{njl} \le 1,
   \quad
   \forall n\in[N],\ l\in\cI_n,
   \nonumber
   \\
   &
   \sum_{l\in\cI_n} h_{nl} \le 1,
   \quad
   \forall n\in[N],
   \nonumber
   \\
   &
   (k_{njl},1,\bgamma_l^\top\ba_{nj}-z_{nl})
   \in\cK_{\exp},
   \quad
   \forall n\in[N],\ l\in\cI_n,\ j\in\cN_l\cap S_n,
   \nonumber
   \\
   &
   (h_{nl},1,q_{nl}-y_n)
   \in\cK_{\exp},
   \quad
   \forall n\in[N],\ l\in\cI_n,
   \nonumber
   \\
   &
   q_{nl} \in \operatorname{MC}_r\big(\lambda_l,z_{nl}\big),
   \quad
   \forall n\in[N],\ l\in\cI_n,
   \nonumber
   \\
   &
   \gamma_{ld} \in \operatorname{MC}_r\big(\rho_l,\beta_d\big),
   \quad
   \forall l\in[L],\ d\in[p],
   \nonumber
   \\
   &
   \rho_l \in
   \big[
        \tfrac{1}{\overline{\lambda}_l^r},
        \tfrac{1}{\underline{\lambda}_l^r}
   \big],
   \quad
   \lambda_l \in
   [\underline{\lambda}_l^r,\overline{\lambda}_l^r],
   \quad
   \forall l\in[L],
   \nonumber
   \\
   &
   \bbt\in\cB,\qquad
   z_{nl}\in[\underline z_{nl}^r,\overline z_{nl}^r],
   \quad \forall n\in[N],\ l\in\cI_n.
   \nonumber
\end{align}
The McCormick notation refers to the envelope
\eqref{eq:mccormick-1}--\eqref{eq:mccormick-4} of the bilinear term \(xy\) over
the node-\(r\) bounds, i.e.,
\([\underline{\lambda}_l^r,\overline{\lambda}_l^r]\times[\underline{z}_{nl}^r,\overline{z}_{nl}^r]\)
for \(q_{nl}\) and
\([\underline{\rho}_l^r,\overline{\rho}_l^r]\times[\underline{\beta}_d,\overline{\beta}_d]\)
for \(\gamma_{ld}\); the cone constraints on \(k_{njl}\) are imposed for all
\(n\in[N]\), \(j\in\cN_l\cap S_n\), and \(l\in\cI_n\). In addition, our
implementation uses additional \emph{perspective cuts}, meaning bounds
based on the convex perspective of log-sum-exp. At exact inclusive values,
$q_{nl}=\lambda_l\log\sum_{j\in\cN_l\cap S_n}e^{v_{nj}/\lambda_l}$.
For fixed utilities $v_{nj}$, this is convex in $\lambda_l$, so tangents give
lower bounds and interval secants give upper bounds. When utilities remain
variable, any imposed cut must be valid jointly over the coefficient and scale
domains; a secant computed only at one fixed utility vector does not by itself
provide such a cut. The validity assumption in
Theorem~\ref{thm:bnb-validity} applies to every added cut.

The equality $\rho_l=1/\lambda_l$ defines the exact lift; the displayed
relaxation retains its reciprocal interval bounds, which need not enforce
that equality at a positive-width node. When the scale interval shrinks to a
point, these reciprocal bounds collapse to the corresponding reciprocal.
Problem \eqref{prob:NL-MC-ECP} is an exponential-cone program because all
nonlinear nonconvex equalities have been replaced by linear McCormick
inequalities, while the remaining nonlinearities are represented by exponential
cones. Since \eqref{prob:NL-MC-ECP} relaxes the feasible region of the original
nested-logit estimation problem, its optimal value gives a valid upper bound on the
maximum log-likelihood over node \(r\).

\paragraph{Lower bound and incumbent update.}

To obtain a feasible solution, we fix \(\blambda\) to a candidate value
\(\hat{\blambda}\) obtained from the relaxation, for example the relaxation
solution or the midpoint of the current box, and solve the original fixed-scale
NL-ECP:
\[
    \max_{\bbt\in\cB}
    \ \cL^{\NL}(\bbt,\hat{\blambda}).
\]
Because \(\hat{\blambda}\) is fixed, this problem is convex and can be solved
globally by an exponential-cone solver. Let $\widehat L_r=\cL^{\NL}(\widehat\bbt_r,\hat\blambda)$ be the
likelihood of the feasible coefficient estimate $\widehat\bbt_r$ obtained by
this solve. It is a valid node-probe lower bound and can improve the global
incumbent.
If its likelihood value improves the incumbent solution, the incumbent is
updated. Improving points are pushed further by a local profile-likelihood
search over \(\blambda\) (re-optimizing \(\bbt\) exactly for each probe)
followed by a joint L-BFGS-B, SLSQP, or trust-constr polishing step.

\paragraph{Branching rule.}

After solving the node relaxation, if the node cannot be pruned, we branch on
one of the dissimilarity parameters. Rather than simply selecting the widest
interval, we pick the nest whose box induces the largest McCormick error in
the bilinear terms \(q_{nl}=\lambda_l z_{nl}\), measured by the score
\begin{equation}\label{eq:branch-score}
    s_l \;=\; \big(\overline{\lambda}_l^r-\underline{\lambda}_l^r\big)
    \sum_{n\in[N]:\,l\in\cI_n} \big(\overline{z}_{nl}^r-\underline{z}_{nl}^r\big),
\end{equation}
Among free scales whose width exceeds $\delta$ or whose endpoint ratio
exceeds $1+\delta$, select a largest-score index $l^\star$. If none is
eligible, the relaxation-error condition and $H\delta+\eta\leq\varepsilon$
ensure that the node passes the objective-gap pruning test after its probe.
Here
\([\underline{z}_{nl}^r,\overline{z}_{nl}^r]\) are the current interval bounds on
\(z_{nl}\). Because the relaxation error of the terms
\(\bbt^\top\ba_{nj}/\lambda_l\) is driven by the ratio
\(\overline{\lambda}_{l}^r/\underline{\lambda}_{l}^r\) rather than by the
width, the interval of \(l^\star\) is split at its geometric mean when this
ratio is large:
\begin{equation}\label{eq:branch-split}
    \lambda_{l^\star}^{\mathrm{split}}
    =
    \begin{cases}
        \sqrt{\underline{\lambda}_{l^\star}^r\,\overline{\lambda}_{l^\star}^r},
        & \text{if }
        \overline{\lambda}_{l^\star}^r/\underline{\lambda}_{l^\star}^r > 4,\\[2pt]
        \tfrac{1}{2}\big(\underline{\lambda}_{l^\star}^r
        +\overline{\lambda}_{l^\star}^r\big),
        & \text{otherwise},
    \end{cases}
\end{equation}
creating the two child nodes
\([\underline{\lambda}_{l^\star}^r,\lambda_{l^\star}^{\mathrm{split}}]\) and
\([\lambda_{l^\star}^{\mathrm{split}},\overline{\lambda}_{l^\star}^r]\).
All other \(\lambda_l\)-intervals are inherited unchanged by the children.
Children are screened by feasibility-based bound tightening (FBBT), which propagates variable bounds through the constraints to shrink a node's box without any solve, and by an analytic interval upper bound $\widetilde{UB}_r$, obtained by
bounding the likelihood over the node intervals, before insertion; at
shallow depths we additionally solve the child relaxation (strong branching),
which both tightens the queue key and provides incumbent probes.

\paragraph{Pruning.}
An \emph{incumbent cut} restricts the node relaxation to objective values at
least $LB$, retaining every point that could improve the incumbent.
Optimization-based bound tightening (OBBT) minimizes and maximizes selected
variables over this relaxation to tighten their interval endpoints.
Let \(UB_r\) denote the upper bound available for node \(r\) --- the optimal
value of the relaxation \eqref{prob:NL-MC-ECP} when it has been solved, and
the analytic interval bound otherwise --- and let \(LB\) denote the value of
the current incumbent. Algorithm~\ref{alg:global-NL} prunes node \(r\) as soon
as any of its bounding steps certifies that the node cannot improve the
incumbent by more than the prescribed tolerance \(\varepsilon>0\):
(i) at selection, if its queue key satisfies \(UB_r\le LB+\varepsilon\);
(ii) if FBBT or reduced-cost bound tightening (RCBT), which uses the dual/reduced-cost information from the relaxation solution to tighten bounds essentially for free, empties the \((\blambda,\bbt,\bgamma)\)-box;
(iii) if the analytic bound satisfies \(\widetilde{UB}_r\le LB+\varepsilon\);
(iv) if the relaxation --- which includes the incumbent cut --- is infeasible,
certifying that no point of the box attains the incumbent value; and
(v) after the lower-bound probe, if \(UB_r\le LB+\varepsilon\) or the internal
gap \(UB_r-\widehat L_r\) falls within the same average-objective tolerance.
Child nodes are subjected to the same tests (FBBT, the analytic bound, and,
at shallow depths, the child relaxation) before insertion, so hopeless nodes
never enter the queue. Since the incumbent value only increases, none of
these tests can discard a point whose log-likelihood exceeds
\(LB+\varepsilon\). When the time limit is reached, $UB$ retains the tolerance-pruning allowance
and every unresolved region's bound, as specified in
Algorithm~\ref{alg:global-NL}. Together with $LB$, it certifies the global
optimality gap for the returned incumbent.
\textcolor{black}{Heuristic searches and polishing steps may be used only to
improve the incumbent. They never justify pruning. A node may be discarded
only by a mathematically valid bound, an infeasibility certificate, or a valid
cut; any unsolved node remains in the global upper-bound accounting.}

\paragraph{Proof of Theorem~\ref{thm:bnb-validity}.}
\label{ec:proof:bnb-validity}
Under the finite bounds on \(\bbt\) and \(\blambda\) (which induce finite
bounds on \(\bz\), \(\bq\), and \(\bgamma\)), every operation of
Algorithm~\ref{alg:global-NL} preserves the invariant that no point with
log-likelihood larger than \(LB+\varepsilon\) is ever discarded: FBBT propagates
valid interval bounds, the perspective and utility cuts are valid inequalities
for the lifted feasible set, and RCBT and OBBT only remove regions whose
relaxation value --- an upper bound on the true likelihood --- provably cannot
exceed the incumbent threshold. Fixing \(\lambda_l=1\) for singleton nests is
also without loss of generality, since \(\lambda_l\) cancels from the
likelihood of a nest containing a single alternative. Although the algorithm
branches on \(\blambda\) alone, this suffices for exactness in the limit:
as the interval of \(\lambda_l\) shrinks in width and in ratio, the interval
of \(\rho_l=1/\lambda_l\) shrinks as well (recall
\(\lambda_l\ge\underline{\lambda}_l>0\)), and the McCormick envelopes of the
bilinear terms \(q_{nl}=\lambda_l z_{nl}\) and \(\gamma_{ld}=\rho_l\beta_d\)
collapse to the exact relations, so the node upper bound converges to the
true maximum of the likelihood over the node box. Because every split removes
a constant fraction of the width or of the log-range of one
\(\lambda_l\)-interval and refinement stops at the branching tolerance, the
search tree is finite. The uniform bound $H\delta+\eta\leq\varepsilon$
ensures that every fully refined box passes the same objective-gap test used
for pruning. Consequently, Algorithm~\ref{alg:global-NL} terminates
after finitely many nodes and returns an incumbent whose log-likelihood lies
within the prescribed tolerance of the global optimum, provided the convex
node relaxations are solved to sufficient accuracy (they are exponential-cone
programs; see the complexity discussion below).

\subsection{Complexity and Practical Implications}
\label{ec:proof:nl-bnb-complexity}
This subsection proves Proposition~\ref{prop:nl-bnb-complexity-main}.

We now discuss the computational complexity of Algorithm~\ref{alg:global-NL}.
The main observation is that all per-node work --- FBBT, the analytic bound,
the node relaxation \eqref{prob:NL-MC-ECP}, RCBT and OBBT, and the
fixed-\(\blambda\) lower-bound probe --- amounts to interval arithmetic plus a
bounded number of exponential-cone programs, each solvable in polynomial time
by a primal--dual interior-point method.

Let
\begin{align*}
    Z &:= \sum_{n=1}^N |S_n|,
    \\
    \Lambda &:= \sum_{n=1}^N
    \left|\{l\in[L]: \cN_l\cap S_n\neq\emptyset\}\right|,
\end{align*}
where \(Z\) is the total number of available alternative appearances and
\(\Lambda\) is the total number of active nest-observation pairs. In the full
availability case, \(Z=NJ\) and \(\Lambda=NL\).

\paragraph{Per-node cost.}
At a node \(r\), the dissimilarity parameters are restricted to the box
\([\underline{\lambda}_l^r,\overline{\lambda}_l^r]\), \(l\in[L]\). The node
relaxation contains \(Z+\Lambda\) exponential-cone blocks (one per constraint
on \(k_{njl}\) and on \(h_{nl}\)), while the McCormick relaxations of
\(q_{nl}=\lambda_l z_{nl}\) and \(\gamma_{ld}=\rho_l\beta_d\) contribute
\(O(\Lambda)\) and \(O(Lp)\) auxiliary variables and linear inequalities,
respectively; the utility and perspective cuts add \(O(Z+\Lambda)\) further
linear rows. Here utility cuts are valid linear inequalities obtained from
the bounds on $v_{nj}=\bbt^\top\ba_{nj}$. The relaxation size is therefore polynomial in \(N\), \(J\),
\(L\), and \(p\), and by the standard complexity theory of interior-point
methods for exponential-cone programming
\citep{NesterovNemirovskii1994,BenTalNemirovski2001} it is solved to accuracy \(\eta>0\)
in
\[
    O\left(
        \sqrt{Z+\Lambda}\log(1/\eta)\cdot
        \operatorname{poly}(N,J,L,p)
    \right)
\]
arithmetic operations. More explicitly, if \(D\) denotes the number of global
variables appearing in the Schur-complement system of the conic relaxation,
then the per-solve cost takes the form
\[
    O\left(
        \sqrt{Z+\Lambda}\log(1/\eta)
        \left(
            Zp^2 + D^3
        \right)
    \right),
\]
where \(D = O\big(p + Lp + L + N + \Lambda\big)\).
The remaining operations of Algorithm~\ref{alg:global-NL} do not change this
order. Each FBBT sweep is \(O(Z+Lp)\) interval arithmetic and a constant
number of sweeps is performed per node; the analytic bound costs
\(O(Z+\Lambda)\); RCBT reuses the duals of the already-solved relaxation at
negligible cost; the fixed-\(\blambda\) lower-bound probe is one additional
(smaller) exponential-cone program; and OBBT, strong branching, and the local
incumbent refinement are capped by fixed time budgets and fire only at the
root, at shallow depths, or at near-prunable nodes. The total work per node is
therefore a constant multiple of the cost of one relaxation solve.

\paragraph{Number of nodes.}
Algorithm~\ref{alg:global-NL} makes nest \(l\) ineligible for further branching once its interval
satisfies both \(\overline{\lambda}_l^r-\underline{\lambda}_l^r\le\delta\) and
\(\overline{\lambda}_l^r/\underline{\lambda}_l^r\le 1+\delta\), where
\(\delta>0\) is a scale-box refinement tolerance. This parameter controls
interval width and ratio; it is distinct from the additive likelihood tolerance
$\varepsilon$ and the experimental relative tolerance
$\varepsilon_{\mathrm{rel}}$. Geometric splits in \eqref{eq:branch-split} partition the
log-range of the interval, and arithmetic bisections partition its width, so
coordinate \(l\) is exhausted after at most
\[
    \left\lceil
    \frac{\ln(\overline{\lambda}_l/\underline{\lambda}_l)}{\ln(1+\delta)}
    \right\rceil
    +
    \left\lceil
    \frac{\overline{\lambda}_l-\underline{\lambda}_l}{\delta}
    \right\rceil
    \;=\;
    O(\kappa_l/\delta)
\]
subintervals, where
\[
    \kappa_l := \max\big\{
        \overline{\lambda}_l-\underline{\lambda}_l,\;
        \ln(\overline{\lambda}_l/\underline{\lambda}_l)
    \big\}.
\]
Therefore, in the worst case, the total number of boxes explored by the
algorithm is bounded by
\[
    O\Big(
        \prod_{l\in\cF}
        \big(1+\lceil \kappa_l/\delta \rceil\big)
    \Big),
\]
which is \textcolor{black}{\(O((1+\kappa/\delta)^{L_f})\)} with \(\kappa:=\max_{l\in[L]}\kappa_l\),
i.e., exponential in the number of \textcolor{black}{free nest scales \(L_f\)}, as expected for spatial
branch-and-bound methods applied to nonconvex optimization problems. Note
that the geometric phase is what keeps \(\kappa_l\) moderate when
\(\underline{\lambda}_l\) is small: controlling the ratio-driven error of the
\(1/\lambda_l\) terms by arithmetic bisection alone would require
\(O\big(1/(\delta\underline{\lambda}_l)\big)\) subintervals.

Combining the node count with the per-node conic complexity, the total
worst-case complexity is
\[
    O\Big(
        \textcolor{black}{(1+\kappa/\delta)^{L_f}}\,
        \sqrt{Z+\Lambda}\,\log(1/\eta)
        \big(
            Zp^2 + D^3
        \big)
    \Big),
\]
with \(D=O(p+Lp+L+N+\Lambda)\). Thus, \textcolor{black}{for fixed \(L_f\), a
fixed compact domain, and fixed branching and conic tolerances}, the algorithm is polynomial in the remaining input
dimensions \(N\), \(J\), and \(p\); conversely, when \(L_f\) is part of the
input, the worst-case complexity is exponential in \(L_f\).

To obtain the approximation-scheme statement, let the requested additive
average-log-likelihood error be $\varepsilon>0$ and choose $\delta(\varepsilon)$ and
$\eta(\varepsilon)$ so that the McCormick relaxation error and conic-solver
error sum to at most $\varepsilon$; the quantitative error condition in the
preceding sample-complexity appendix gives one sufficient choice. For fixed
$L_f$ and a fixed compact domain, the factor
$(1+\kappa/\delta(\varepsilon))^{L_f}$ has no additional dependence on
$N$ through objective normalization. Any dependence of $H$ on the feature or
alternative dimensions must still be included; the stated polynomial-time
conclusion requires it to be polynomial when $L_f$ is fixed. Hence, for every fixed $\varepsilon$, Algorithm~\ref{alg:global-NL}
runs in time polynomial in $N$, $J$, and $p$ and returns an estimate whose
log-likelihood is within $\varepsilon$ of the global maximum. This is an
additive PTAS. It is not claimed to be strongly polynomial, and it is not a
multiplicative scheme because the maximized log-likelihood is generally
negative. If $L_f$ is part of the input, the displayed exponential dependence
remains.

This fixed-parameter structure is important in practice.  The number of nests
is a behavioral modeling choice and is usually small (often two or three),
whereas the number of observations and available alternatives can be orders of
magnitude larger.  Increasing the dataset size therefore enlarges the convex
ECP solved at each node but does not add a new dimension to the combinatorial
search.  Adding a nest is qualitatively different: it adds a branching
coordinate and can multiply the worst-case number of boxes by
$O(\kappa/\delta)$.  Practitioners should consequently use the smallest nesting
structure supported by behavioral or statistical evidence, fix singleton nests
at $\lambda_l=1$, and exploit parallelism and sparsity as $N$ grows.  The bound
also explains the empirical pattern in the appendix: runtime per node grows
with sample size, while the number of nodes is driven primarily by the nesting
structure and relaxation tightness.

This behavior is consistent with the division of labor in
Algorithm~\ref{alg:global-NL}: the exponential-cone solver provides
polynomial-time bounds for every convex subproblem, while the global search
over the nonconvex dissimilarity parameters requires, in the worst case, an
exponential number of nodes. The acceleration devices do not affect the
worst-case order --- FBBT and RCBT are low-order polynomial per node, while
OBBT, strong branching, and incumbent refinement are time-budgeted --- but in
practice they shrink the boxes, the perspective and utility cuts tighten every
relaxation, and the warm start and incumbent polishing raise \(LB\) early, so
the number of explored nodes typically remains far below the worst-case
bound. Finally, since every pruned or open node carries a valid upper bound,
Algorithm~\ref{alg:global-NL} is an \emph{anytime} method: upon interruption
it returns the incumbent together with the absolute gap $UB-LB$ and,
when $UB\neq0$, the relative gap $g_{\mathrm{rel}}=(UB-LB)/|UB|$.

\section*{Appendix: Experimental Details and Additional Results}
The benchmark definitions and statistical abbreviations are given here before
Table~\ref{tab:nest-attr-evidence}. LR denotes the likelihood-ratio statistic,
$p_{\mathrm{LR}}$ its p-value, and LM the Lagrange-multiplier (score) test.
BIC is the Bayesian information criterion; $\Delta\mathrm{BIC}$ is BIC for
MNL minus BIC for NL, so a positive value favors NL. A Wald statistic is a
coefficient estimate divided by its estimated standard error; $|t|$ is its
absolute value. In the boundary-mixture notation, $\chi^2_k$ denotes the
chi-squared distribution with $k$ degrees of freedom, with $\chi^2_0$ a point
mass at zero.

Heating alternatives GC, GR, EC, ER, and HP mean gas central, gas room,
electric central, electric room, and heat pump. PT is public transport; SM
means slow modes in Optima and the Swissmetro alternative in SwissMetro.
Yogurt brands Weight, Yop., Hil., and Dan. mean Weight Watchers, Yoplait,
Hiland, and Dannon. Cracker labels Sun., Priv., Keeb., and Nab. mean Sunshine,
private label, Keebler, and Nabisco. Attribute names COST, TIME, PRICE, CATCH,
and DIST mean monetary cost, travel time, price, catch rate, and travel
distance; IC and OC are installation and annual operating costs; IVT, OVT,
and FREQ are in-vehicle time, out-of-vehicle time, and service frequency;
DISP and FEAT indicate a product display and feature advertisement.
The vertical bar in a nest structure separates nests; braces enclose the
alternatives in each nest.

Table~\ref{tab:nest-attr-evidence} reports the nest-selection evidence and retained attributes. Nest selection is repeated
independently within each training subsample; no held-out choice is used to
select a partition, attributes, or tuning parameters. We enumerate partitions into two or three nonempty nests with at least one
nonsingleton nest ($111$ candidate structures in total). The all-singleton
partition is excluded for each three-alternative dataset, Optima and
SwissMetro, because it is identically MNL and has no free dissimilarity.
These two exclusions reduce the unrestricted count from $113$ to $111$;
 a structure is \emph{admissible} if the score (LM) test at
$\lambda{=}1$ points in the nesting-consistent direction and every fitted free
dissimilarity parameter is interior, and among admissible structures the one
with the smallest boundary-corrected likelihood-ratio $p$-value against MNL
is selected (when LR statistics differ by less than $2$, ties are resolved
toward fewer nests).
For summed fitted training log-likelihoods $\mathcal L^{\mathrm{fit}}_{\mathrm{NL}}$
and $\mathcal L^{\mathrm{fit}}_{\mathrm{MNL}}$, the statistic is
$\mathrm{LR}=2(\mathcal L^{\mathrm{fit}}_{\mathrm{NL}}-
\mathcal L^{\mathrm{fit}}_{\mathrm{MNL}})$ against the same-specification MNL model; $p_{\mathrm{LR}}$ is its boundary-corrected p-value ($\tfrac12\chi^2_0+\tfrac12\chi^2_1$ per
dissimilarity parameter); $\Delta$BIC $=\mathrm{LR}-L_f\ln N$ for $L_f$ free
dissimilarity parameters ($>0$ favors the nested model). The NL\,/\,MNL column
counts, among all structures scanned for that dataset, those whose
fitted estimate is genuinely nested versus those that collapse to MNL (anti-nested score direction); the final row totals these counts: $35$ of $111$ candidate structures require nested logit, while $76$ collapse exactly to
MNL. 

Every attribute is statistically significant (Wald $|t|$ shown, all well above the $5\%$ threshold of $1.96$), so all are retained and none are dropped.

\begin{table*}[tbp]
\centering
\footnotesize
\setlength{\tabcolsep}{4.5pt}
\resizebox{\textwidth}{!}{%
\begin{tabular}{@{}llrrrc l@{}}
\toprule
Dataset ($N_{\mathrm{full}}$) & Selected nest structure & LR & $p_{\mathrm{LR}}$ & $\Delta$BIC & NL\,/\,MNL & Attributes ($|t|$) \\
\midrule
Mode (453) & $\{$Car, Carpool$\}\,|\,\{$Bus$\}\,|\,\{$Rail$\}$ & 28.02 & $6.0\times10^{-8}$ & $+21.9$ & 4\,/\,9
& COST (4.0), TIME (14.3) \\
Heating (900) & $\{$GC, GR, EC, HP$\}\,|\,\{$ER$\}$ & 44.66 & $1.2\times10^{-11}$ & $+37.9$ & 9\,/\,31
& IC (17.7), OC (14.2) \\
Fishing (1{,}182) & $\{$Charter$\}\,|\,\{$Beach, Pier, Boat$\}$ & 153.6 & $1.4\times10^{-35}$ & $+146.6$ & 5\,/\,8
& PRICE (16.7), CATCH (10.7) \\
Optima (1{,}906) & $\{$Car$\}\,|\,\{$PT, SM$\}$ & 10.69 & $5.4\times10^{-4}$ & $+3.1$ & 1\,/\,2
& TIME (9.8), COST (10.1), DIST (16.8) \\
Yogurt (2{,}412) & $\{$Weight, Yop., Hil.$\}\,|\,\{$Dan.$\}$ & 378.4 & $1.4\times10^{-84}$ & $+370.6$ & 5\,/\,8
& PRICE (11.0), FEAT (9.0) \\
ModeCanada (2{,}779) & $\{$Train$\}\,|\,\{$Air, Bus, Car$\}$ & 882.8 & $2.7\times10^{-194}$ & $+874.9$ & 4\,/\,9
& COST (8.5), IVT (23.1), OVT (35.6), FREQ (7.8) \\
Cracker (3{,}292) & $\{$Sun., Priv.$\}\,|\,\{$Keeb.$\}\,|\,\{$Nab.$\}$ & 484.9 & $9.4\times10^{-108}$ & $+476.8$ & 5\,/\,8
& PRICE (10.9), DISP (13.8), FEAT (7.1) \\
SwissMetro (5{,}607) & $\{$Train, SM$\}\,|\,\{$Car$\}$ & 246.3 & $8.5\times10^{-56}$ & $+237.6$ & 2\,/\,1
& TIME (45.5), COST (23.3) \\
\midrule
All (111 structures) & & & & & \textbf{35}\,/\,\textbf{76} & \\
\bottomrule
\end{tabular}}
\caption{Selected nest structures, model-comparison statistics, and attribute significance for the eight benchmarks.}
\label{tab:nest-attr-evidence}
\end{table*}
All models are estimated by maximum likelihood under a common $1$\,hour
wall-clock limit. The convex MNL likelihood is solved by MNL-ECP with
\textsc{Mosek} to a relative/absolute
optimality gap of $10^{-4}$; being convex, every instance is certified to the
stated global tolerance (final gap $<10^{-5}$) within seconds, well inside the limit.
The nonconvex nested-logit likelihood is solved by NL-B\&B to a
\emph{certified} relative optimality gap of $10^{-2}$: the search
maintains a global lower bound (a feasible incumbent) and an upper bound (the
relaxation value), and terminates when $(UB-LB)/|UB|
\le 10^{-2}$ or the time limit is reached. Each node's exponential-cone
relaxation is solved to a relative gap of $10^{-4}$. L-BFGS-B, SLSQP, and trust-constr each use five random
restarts and solver tolerance $10^{-4}$ as the local-solver baselines. All 240 real-data subsample instances reach the prescribed
$\varepsilon_{\mathrm{rel}}=10^{-2}$ certificate within the nominal one-hour
budget. Their certified incumbents provide the reference values for the local
solver comparisons; no extended-run reference is used. Recorded elapsed times
can exceed the nominal budget because of system delays. For presentation,
training-fraction mean runtimes exceeding the budget are displayed as
$3{,}600$ seconds; the main table averages these displayed fraction means.
This convention caps reported aggregate times rather than recomputing
per-instance measurements, and leaves objective values, gaps, and node counts
unchanged. The separate full-sample convergence
runs use the tighter $\varepsilon_{\mathrm{rel}}=10^{-3}$ tolerance and may
stop before attaining it; these runs are not among the 240 subsample instances.
\paragraph{Data:} We use eight public datasets to evaluate both MNL and
nested-logit. Cracker, Fishing, Heating, Mode, and Yogurt come from the Ecdat
collection (\url{https://vincentarelbundock.github.io/Rdatasets/articles/data.html}).
Optima and SwissMetro are obtained from TRANSP-OR
(\url{https://github.com/mncosta/biogeme_tutorial/raw/master/Data/optima.dat} and
\url{https://github.com/timothyb0912/pylogit/raw/master/examples/data/swissmetro.dat}),
and ModeCanada is obtained from the \texttt{mlogit} R package
(\url{https://github.com/gsbDBI/torch-choice/raw/main/tutorials/public_datasets/ModeCanada.csv}).

A dataset-specific script downloads each source and converts it to a common
conditional-logit representation. We retain observations with a reported choice
and a complete menu of available alternatives, and discard records that do not
define a valid choice situation. For example, we restrict ModeCanada to cases
with all four alternatives and reshape its long-format records to one row per
traveller. For SwissMetro, we retain commuter and business trips for which all
three modes are available. Alternatives are represented by their observed
attributes, such as travel time, cost, price, and other level-of-service or
product characteristics. Monetary and time variables are rescaled by a fixed
constant, typically $1/100$, to bring their magnitudes close to one. This
transformation is applied before sampling and is shared by the full dataset, all
subsamples, and both logit specifications.

From each cleaned dataset, we draw 10 independent random subsamples at each of
three sampling rates ($20\%$, $50\%$, and $80\%$), using a fixed seed for every
replication. This procedure yields 30 estimation instances per dataset. For each
instance, the observations not selected for estimation form the held-out set.
We also retain the complete cleaned dataset as the common sampling pool. All
instances use the same plain-text schema: a header specifying the numbers of
individuals, alternatives, attributes, and nests; the attribute matrix; each
individual's availability set and chosen alternative; and the nest partition.
Both NL-B\&B and MNL-ECP read these files without further transformation.

\subsection*{Additional Results}
Table~\ref{tab:nl-cert} gives the detailed certification results.  Datasets are columns, and both blocks are grouped by sub-sample fraction
($0.2N_{\mathrm{full}}/0.5N_{\mathrm{full}}/0.8N_{\mathrm{full}}$, ten instances each). \emph{Top:} NL-B\&B certifies every instance to the prescribed $1\%$ relative tolerance; per fraction we report the number of B\&B nodes and recorded wall-clock seconds, with displayed
times capped at the nominal $3{,}600$-second budget when system delays cause
overruns. \emph{Bottom:} each
local solver is run from five restarts per instance; the best, median and worst
restart objectives give per-instance percentage shortfalls computed as
$100(\mathcal L_{\mathrm{ref}}-\mathcal L_{\mathrm{loc}})/|\mathcal L_{\mathrm{ref}}|$,
where $\mathcal L_{\mathrm{ref}}$ is the certified training reference and
$\mathcal L_{\mathrm{loc}}$ the corresponding restart objective ---
the outcome of keeping the best run, of a typical run, and of an unlucky one.
Per fraction, each cell is the mean of that gap over the ten instances with the
single-worst instance in parentheses, e.g.\ ``0.55 (4.31)''; ``--'' marks a gap
below $0.01\%$. Multi-starting cannot repair the failure: on a Heating
instance the best of five restarts still lands $4.31\%$ below the reference for
all three solvers, while an individual run strays as far as $7.24\%$ (L-BFGS-B),
$22.10\%$ (trust-constr) and $98.84\%$ (SLSQP) on some ModeCanada instances.
The trap is invisible without the certificate, and which solver is safe is not
knowable a priori.

\begin{table*}[tbp]
\centering
\resizebox{\textwidth}{!}{%
\begin{tabular}{@{}ll rrrrrrrr@{}}
\toprule
 & & \multicolumn{1}{c}{Mode} & \multicolumn{1}{c}{Heating} & \multicolumn{1}{c}{Fishing} & \multicolumn{1}{c}{Optima} & \multicolumn{1}{c}{Yogurt} & \multicolumn{1}{c}{ModeCanada} & \multicolumn{1}{c}{Cracker} & \multicolumn{1}{c}{SwissMetro} \\
\midrule
\multicolumn{10}{@{}l}{\textit{NL-B\&B --- \textcolor{black}{certified to the prescribed global tolerance}; mean search cost per sub-sample}}\\
\midrule
\multirow{2}{*}{$0.2N_{\mathrm{full}}$} %& $N$ & 91 & 180 & 236 & 381 & 482 & 556 & 658 & 1{,}121 \\
 & Nodes & 181 & 261 & 38 & 629 & 30 & 142 & 40 & 204 \\
 & Time (s) & 41.82 & 81.81 & 26.90 & 512.51 & 37.36 & 284.16 & 98.07 & 858.00 \\
\addlinespace
\multirow{2}{*}{$0.5N_{\mathrm{full}}$} %& $N$ & 226 & 450 & 591 & 953 & 1{,}206 & 1{,}390 & 1{,}646 & 2{,}804 \\
 & Nodes & 219 & 294 & 31 & 487 & 24 & 123 & 62 & 180 \\
 & Time (s) & 139.00 & 234.03 & 58.33 & 958.87 & 90.50 & 826.30 & 467.49 & 2{,}182.02 \\
\addlinespace
\multirow{2}{*}{$0.8N_{\mathrm{full}}$} %& $N$ & 362 & 720 & 946 & 1{,}525 & 1{,}930 & 2{,}223 & 2{,}634 & 4{,}486 \\
 & Nodes & 194 & 275 & 44 & 550 & 24 & 173 & 74 & 267 \\
 & Time (s) & 210.53 & 380.47 & 124.49 & 1{,}888.56 & 153.60 & 1{,}673.87 & 723.00 & 3{,}600.00 \\
\midrule
\multicolumn{10}{@{}l}{\textit{L-BFGS-B --- gap below the certified reference per sub-sample}}\\
\midrule
\multirow{3}{*}{$0.2N_{\mathrm{full}}$} & Best restart gap (\%) & -- & 0.74 (4.31) & -- & -- & 0.15 (0.45) & -- & 0.02 (0.11) & -- \\
 & Median restart gap (\%) & 0.19 (1.38) & 1.76 (4.54) & 0.01 (0.03) & 0.35 (1.17) & 0.83 (1.54) & 0.02 (0.05) & 0.13 (0.38) & -- \\
 & Worst restart gap (\%) & 0.57 (3.33) & 2.82 (5.69) & 0.02 (0.14) & 0.98 (2.86) & 2.17 (4.28) & 0.08 (0.25) & 0.29 (0.94) & 0.01 (0.04) \\
\addlinespace
\multirow{3}{*}{$0.5N_{\mathrm{full}}$} & Best restart gap (\%) & -- & 0.58 (2.09) & -- & -- & 0.14 (0.28) & -- & 0.02 (0.05) & -- \\
 & Median restart gap (\%) & 0.03 (0.13) & 1.41 (2.10) & -- & 0.14 (0.31) & 0.62 (1.15) & 0.18 (1.54) & 0.17 (0.43) & 0.00 (0.01) \\
 & Worst restart gap (\%) & 0.14 (0.66) & 2.34 (3.48) & 0.01 (0.02) & 0.47 (1.12) & 1.54 (4.05) & 0.86 (7.24) & 0.43 (1.78) & 0.01 (0.05) \\
\addlinespace
\multirow{3}{*}{$0.8N_{\mathrm{full}}$} & Best restart gap (\%) & -- & 0.68 (1.57) & -- & -- & 0.13 (0.45) & -- & 0.02 (0.07) & -- \\
 & Median restart gap (\%) & -- & 1.59 (1.82) & -- & 0.11 (0.21) & 0.64 (2.00) & 0.02 (0.05) & 0.12 (0.24) & -- \\
 & Worst restart gap (\%) & 0.01 (0.01) & 2.21 (2.80) & 0.01 (0.05) & 0.46 (0.85) & 1.60 (5.30) & 0.07 (0.26) & 0.33 (0.45) & 0.02 (0.04) \\
\midrule
\multicolumn{10}{@{}l}{\textit{SLSQP --- gap below the certified reference per sub-sample}}\\
\midrule
\multirow{3}{*}{$0.2N_{\mathrm{full}}$} & Best restart gap (\%) & -- & 0.55 (4.31) & -- & -- & -- & -- & -- & -- \\
 & Median restart gap (\%) & -- & 0.62 (4.31) & -- & -- & 0.04 (0.38) & 22.85 (79.70) & -- & -- \\
 & Worst restart gap (\%) & -- & 0.89 (4.31) & -- & -- & 0.19 (1.84) & 28.23 (95.15) & -- & -- \\
\addlinespace
\multirow{3}{*}{$0.5N_{\mathrm{full}}$} & Best restart gap (\%) & -- & 0.45 (2.09) & -- & -- & -- & -- & -- & -- \\
 & Median restart gap (\%) & -- & 0.53 (2.09) & -- & -- & -- & 46.12 (78.87) & -- & -- \\
 & Worst restart gap (\%) & -- & 0.75 (2.32) & -- & -- & 0.01 (0.01) & 56.59 (94.91) & -- & -- \\
\addlinespace
\multirow{3}{*}{$0.8N_{\mathrm{full}}$} & Best restart gap (\%) & -- & 0.16 (1.57) & -- & -- & -- & -- & -- & -- \\
 & Median restart gap (\%) & -- & 0.19 (1.57) & -- & -- & -- & 53.97 (95.20) & -- & -- \\
 & Worst restart gap (\%) & -- & 0.27 (1.57) & -- & -- & 0.01 (0.02) & 58.33 (98.84) & -- & -- \\
\midrule
\multicolumn{10}{@{}l}{\textit{trust-constr --- gap below the certified reference per sub-sample}}\\
\midrule
\multirow{3}{*}{$0.2N_{\mathrm{full}}$} & Best restart gap (\%) & -- & 0.43 (4.31) & -- & -- & -- & -- & -- & -- \\
 & Median restart gap (\%) & -- & 0.66 (4.32) & 0.35 (1.16) & -- & -- & 0.54 (5.37) & -- & 0.05 (0.45) \\
 & Worst restart gap (\%) & -- & 1.20 (4.32) & 1.67 (5.53) & -- & 0.01 (0.02) & 2.21 (22.10) & -- & 0.22 (2.22) \\
\addlinespace
\multirow{3}{*}{$0.5N_{\mathrm{full}}$} & Best restart gap (\%) & -- & 0.08 (0.83) & -- & -- & -- & -- & -- & -- \\
 & Median restart gap (\%) & 0.14 (0.79) & 0.39 (1.68) & -- & -- & -- & 1.61 (7.64) & -- & 0.06 (0.59) \\
 & Worst restart gap (\%) & 0.67 (3.82) & 0.48 (2.09) & -- & -- & -- & 5.34 (18.47) & -- & 0.29 (2.88) \\
\addlinespace
\multirow{3}{*}{$0.8N_{\mathrm{full}}$} & Best restart gap (\%) & -- & -- & -- & -- & -- & -- & -- & -- \\
 & Median restart gap (\%) & 0.07 (0.71) & 0.13 (1.31) & 0.12 (1.18) & 0.00 (0.04) & -- & 2.07 (8.25) & -- & 0.06 (0.58) \\
 & Worst restart gap (\%) & 0.34 (3.45) & 0.18 (1.79) & 0.57 (5.65) & 0.02 (0.19) & -- & 7.16 (18.36) & -- & 0.28 (2.84) \\
\bottomrule
\end{tabular}%
}
\caption{NL-B\&B search cost and local-solver shortfalls by dataset and training fraction. Entries are means, with worst-case shortfalls in parentheses.}
\label{tab:nl-cert}
\end{table*}

\begin{figure*}[t]
    \centering
    \includegraphics[width=\linewidth]{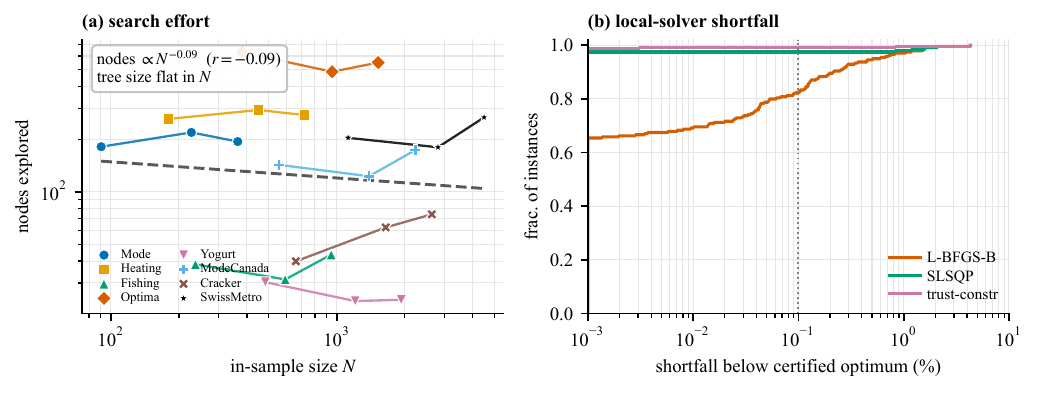}
    \caption{Search effort versus training size (left) and cumulative distribution of local-solver shortfalls across 240 instances (right).}
    \label{fig:bnb}
\end{figure*}

Figure~\ref{fig:bnb} relates search cost to the benefit of certification.
The left panel plots B\&B nodes explored against in-sample size $N\in[90,4484]$; each point averages ten random sub-samples of one dataset
(colours in the legend), the dashed line is the pooled power-law fit and $r_{\mathrm{corr}}$ denotes
the reported correlation coefficient (printed as $r$ inside the plot). The tree
size is essentially \emph{independent} of $N$ ($\propto N^{-0.09}$, $r_{\mathrm{corr}}=-0.09$):
it is set by the nonconvexity of the likelihood, which differs across datasets
(Optima explores ${\sim}600$ nodes, Yogurt and Fishing fewer than $45$), not by
how many individuals are sampled. Only the per-node cost grows with $N$
($\propto N^{1.18}$, $r_{\mathrm{corr}}=0.98$, as each node solves an exponential-cone
relaxation with $O(N)$ conic constraints), giving near-linear total scaling
(${\approx}N^{1.09}$): $N$ enters through the tractable relaxation, never
through the combinatorial search. 

The right panel shows the empirical cumulative distribution function, over
all $240$ instances, of the shortfall below the certified
reference incurred by the \emph{best} of five restarts of each local
solver. Multi-start L-BFGS-B matches the reference within the $0.01\%$ reporting threshold on only
$65\%$ of instances and misses it by more than $0.1\%$ on $18\%$; SLSQP and
trust-constr are more reliable ($97.5\%$ and $98.8\%$ within the same reporting threshold) but all three fall below the reference
on the same worst instance, by $4.31\%$. Every failing run reports normal
convergence, so without the certificate the error is invisible.

\begin{figure*}[t]
    \centering
    \includegraphics[width=\linewidth]{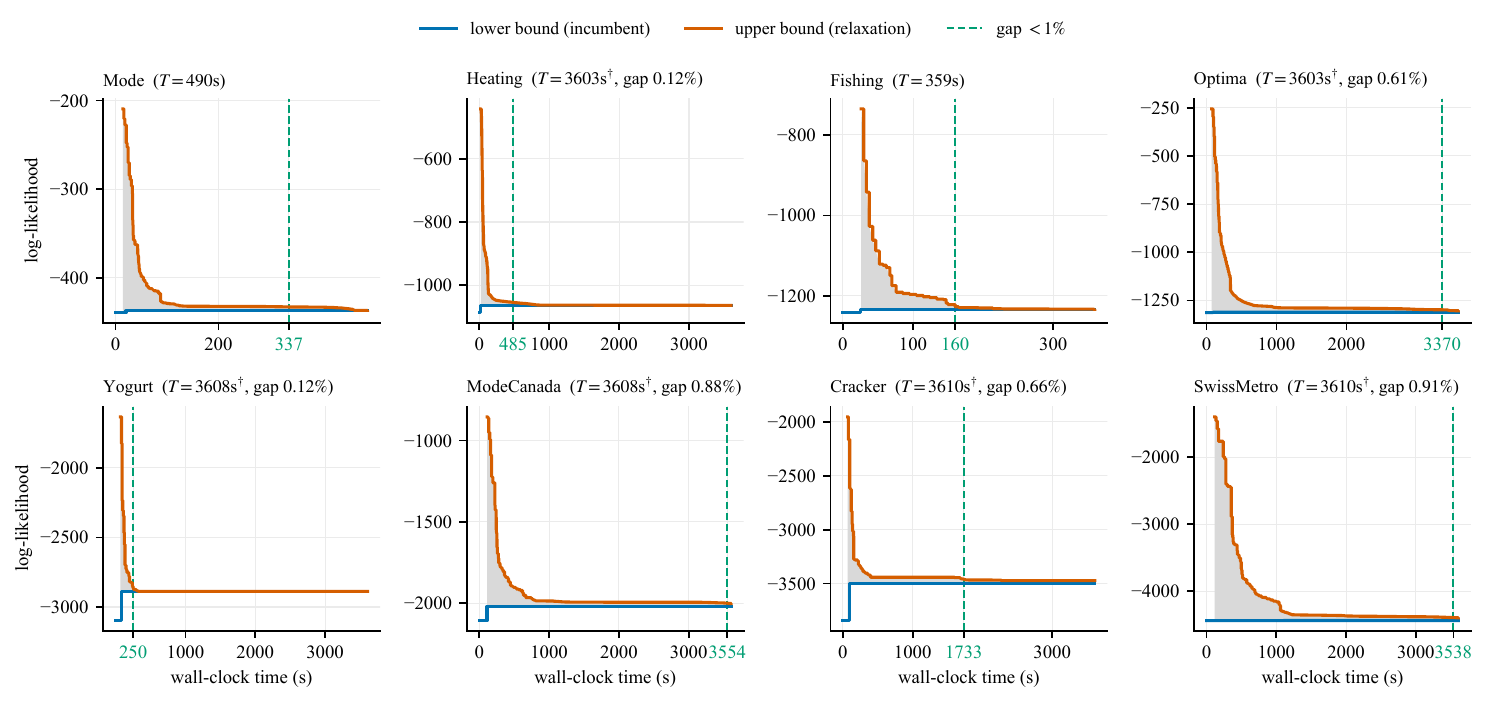}
    \caption{Anytime likelihood bounds on the eight full datasets. Shading marks the certified gap; green lines mark the first $1\%$ certificate.}
    \label{fig:convergence-panels}
\end{figure*}

Figure~\ref{fig:convergence-panels} shows the full-sample bound trajectories,
with one panel per dataset ($N=N_{\mathrm{full}}=453$--$5{,}607$) and target
$\varepsilon_{\mathrm{rel}}=10^{-3}$. In each panel, the blue step curve is the certified lower bound (the incumbent
log-likelihood, non-decreasing), the orange step curve the certified upper
bound (including open-node bounds and the tolerance-pruning allowance), and the
shaded band between them the certified optimality gap; any point outside the
band is provably not the global optimum. The green dashed line, with its time
marked in green on the horizontal axis, is the first moment the gap reaches
$1\%$---the tolerance targeted throughout the rest of the paper---so from that
line onward the incumbent is guaranteed within $1\%$ of the global maximum.
Panel titles give the recorded elapsed time $T$ in seconds, including any
system-delay overrun; these trace labels are uncapped. A dagger marks the six runs stopped
by the one-hour budget, for which the gap attained is also reported (runs that
terminate on their own, Fishing and Mode, are certified to $\varepsilon_{\mathrm{rel}}$ by
construction and need no such number).
Two features recur across all eight datasets. First, the blue curve is flat
almost from the start: the incumbent reaches its final value after $1$--$19\%$
of the runtime (median $4\%$; Heating after $20$\,s of $3{,}603$\,s, Cracker
after $97$\,s of $3{,}610$\,s), and in seven of eight cases the very first
improvement over the warm start is the last. The estimate is therefore
available almost immediately; everything after is proof, not search. Second,
all of the visible progress is in the orange curve, which must travel a long
way---from $-1{,}954$ to $-3{,}473$ on Cracker, from $-1{,}397$ to $-4{,}396$
on SwissMetro---before the band closes. The width of that band is
what distinguishes a certified estimate from a merely converged one: a local
solver returns the blue curve alone, with no way to tell how much of the shaded
region remains.

\end{document}